\documentclass[twoside]{article}
\usepackage{graphicx,amssymb,mathrsfs,amsmath}
\usepackage[numbers,sort&compress]{natbib}
\catcode`@=11
\long\def\@makefntext#1{\noindent #1}
\newskip\tabcentering \tabcentering=1000pt plus 1000pt minus 1000pt
\def\MCH#1#2{\setbox0=\hbox{\raise#1\hbox{#2}}\smash{\box0}}

\def\CR{\cr\noalign{\vspace{1mm} \hrule \vspace{1mm}}}
\def\@evenfoot{}\def\@oddfoot{}

\def\@evenhead{\hbox to\textwidth{\footnotesize\rm\thepage \hfill
{\it Tiexin Guo, Shien Zhao, Xiaolin Zeng}}} 

\def\@oddhead{\hbox to \textwidth{\footnotesize{\it
Random convex analysis (I): separation and Fenchel-Moreau duality in random locally convex modules } \hfill\thepage}}

\def\sec#1{\vskip 3mm\leftline{\bf #1}\vskip 1mm}
\def\subsec#1{\vskip 2mm\leftline{#1}\vskip 1mm}
\def\th#1{\vskip 1mm\noindent{\bf #1}\quad}

\renewcommand{\topfraction}{1}
\renewcommand{\bottomfraction}{1}
\renewcommand{\textfraction}{0}
\renewcommand{\floatpagefraction}{0}
\catcode`@=12

\def\R{{\Bbb R}}

\def\bc{\begin{center}}
\def\ec{\end{center}}
\def\no{\noindent}
\def\hang{\hangindent\parindent}
\def\textindent#1{\indent\llap{\qquad #1\ \ \enspace}\ignorespaces}
\def\ref{\par\hang\textindent}

\def\dl{\displaystyle\lim}

\begin{document}
\abovedisplayskip=6pt plus 1pt minus 1pt \belowdisplayskip=6pt
plus 1pt minus 1pt
\thispagestyle{empty} \vspace*{-1.0truecm} \noindent
\vskip 10mm

\bc{\Large\bf Random convex analysis (I): separation and Fenchel-Moreau duality in random locally convex modules

\footnotetext{\footnotesize $^{*}$Corresponding author}} \ec

\vskip 5mm
\bc{\bf Tiexin Guo$^{1,*}$, Shien Zhao$^{2}$, Xiaolin Zeng$^{3}$}\\
{\small $^{1}$School of Mathematics and Statistics, Central South University, Changsha {\rm 410083}, China;\\
$^{2}$Elementary Educational College, Capital Normal University, Beijing {\rm 100048}, China;\\
$^{3}$School of Mathematics and Statistics, Chongqing Technology and Business University, Chongqing {\rm 400067}, China\\
\footnotesize
E-mail: tiexinguo@csu.edu.cn; zsefh@cnu.edu.cn; xlinzeng@163.com}\ec

\def\sec#1{\vspace{2mm}\noindent{{\bf #1}}\vspace{0.5mm}}
\def\subsec#1{\vspace{2mm}\leftline{\bf #1}} 
\def\th#1{\vspace{1mm}\noindent{\bf #1}\quad } 
\def\pf#1{\vspace{1mm}\noindent{\it #1}\quad}
\renewcommand{\topfraction}{1} \renewcommand{\bottomfraction}{1}
\renewcommand{\textfraction}{0} \renewcommand{\floatpagefraction}{0}
\floatsep=0pt \textfloatsep=0pt \intextsep=0pt \catcode`@=12
\def\leq{\leqslant}
\def\geq{\geqslant}
\def\R{{\Bbb R}}  \def\N{{\Bbb N}}  \def\Q{{\Bbb Q}}  \def\O{{\Bbb O}} \def\Z{{\Bbb Z}}
\def\C{{\Bbb C}}  \def\hml{\end{document}}  \newsymbol\wjzhml 203F \def\no{\noindent}
\def\CR{\cr\noalign{\vspace{1mm} \hrule \vspace{1mm}}}

\abovedisplayskip=3pt plus 1pt minus 1pt 
\belowdisplayskip=3pt plus 1pt minus 1pt 

\def\le{\leqslant}
\def\ge{\geqslant}
\def\dl{\displaystyle}

\newtheorem{theorem}{Theorem}[section]%
\newtheorem{corollary}[theorem]{Corollary}%
\newtheorem{lemma}[theorem]{Lemma}%
\newtheorem{proposition}[theorem]{Proposition}%
\newtheorem{problem}[theorem]{Problem}
\newtheorem{definition}[theorem]{Definition}%
\newtheorem{remark}[theorem]{Remark}%
\newtheorem{example}[theorem]{Example}
\newtheorem{claim}[theorem]{Claim}



\vspace{8true mm}

\renewcommand{\baselinestretch}{1.9}\baselineskip 19pt

\baselineskip 12pt \renewcommand{\baselinestretch}{1.18}
\noindent{{\bf Abstract}\small\hspace{2.8mm} 
To provide a solid analytic foundation for the module approach to conditional risk measures, our purpose is to establish a complete random convex analysis over random locally convex modules by simultaneously considering the two kinds of topologies (namely the $(\varepsilon,\lambda)$--topology and the locally $L^0$-- convex topology). This paper is focused on the part of separation and Fenchel-Moreau duality in random locally convex modules. The key point of this paper is to give the precise relation between random conjugate spaces of a random locally convex module under the two kinds of topologies, which enables us to not only give a thorough treatment of separation between a point and a closed $L^{0}$-convex subset but also establish the complete Fenchel-Moreau duality theorems in random locally convex modules under the two kinds of topologies.\vspace{0mm} }

\vspace{1mm} \no{\footnotesize{\bf Keywords:Random locally convex module,
$(\varepsilon,\lambda)$--topology, Locally $L^{0}$--convex topology,
Random conjugate spaces, Separation, Lower semicontinuous $L^{0}$--convex function, Fenchel-Moreau duality

}}

\no{\footnotesize{\bf MSC(2000):\hspace{2mm}46A20, 46A22, 46A55, 46H25.}
 \vspace{2mm}
\baselineskip 15pt
\renewcommand{\baselinestretch}{1.22}
\parindent=10.8pt  
\rm\normalsize\rm

\section{Introduction}

Random metric theory (or, also aptly called random functional analysis) is based on the idea of randomizing the classical space theory of functional analysis. Random normed modules (briefly, $RN$ modules) and random locally convex modules (briefly, $RLC$ modules) together with their random conjugate spaces were naturally formed and have been the central theme in the course of the development of random metric theory, cf. \cite{TXG-Master,TXG-PHD,TXG-Extension,TXG-Module,TXG-Radon,TXG-basic,TXG-Sur,TXG-Alao,TXG-JFA,Guotx-onsome,TXG-SBL}. Classical convex analysis (e.g, see \cite{ET}) is the analytic foundation for convex risk measures, cf. \cite{ADEH,Delbaen,FS,Follmer-S,Fritt-R}. However, it is no longer universally applicable to $L^{0}$--convex (or conditional convex) conditional risk measures (in particular, those defined on the model spaces of unbounded financial positions). Just to overcome the obstacle, D. Filipovi\'{c}, et.al presented the module approach to conditional risk, cf, \cite{FKV,FKV-appro}. Let $(\Omega,\mathcal{E},P)$ be a probability space, $\mathcal{F}$ a sub--$\sigma$--algebra of $\mathcal{E}$, $L^{0}(\mathcal{F})$ ($\bar{L}^{0}(\mathcal{F})$) the set of real (extended real)-valued $\mathcal{F}$-measurable random variables on $\Omega$, $L^{p}(\mathcal{E})$ $(1\leq p\leq+\infty)$ the classical function space and $L^{p}_{\mathcal{F}}(\mathcal{E})$ the  $L^{0}(\mathcal{F})$-module generated by $L^{p}(\mathcal{E})$, which can be made into a random normed module in a natural way, see Example 2.3 of this paper for the construction of the random normed module $L^{p}_{\mathcal{F}}(\mathcal{E})$. The so-called module approach is to choose $L^{p}_{\mathcal{F}}(\mathcal{E})$ as the model space, namely define an $L^{0}(\mathcal{F})$-convex conditional risk measure to be a proper $L^{0}(\mathcal{F})$-convex cash-invariant and monotone function from $L^{p}_{\mathcal{F}}(\mathcal{E})$ to $\bar{L}^{0}(\mathcal{F})$ and further develop conditional risk measures under the definition. Analysis of the conditional risk measures naturally goes beyond the scope of the classical convex analysis of a proper lower semicontinuous extended real-valued convex function on a locally convex space. This requires a complete random convex analysis, namely a module analogue of classical convex analysis in random metric theory, in order to deal with the Fenchel-Moreau dual representation, continuity and subdifferentiability for such a lower semicontinuous $L^{0}(\mathcal{F})$-convex conditional risk measure. Since classical convex analysis is of wide use in various disciplines of mathematics, the study of random convex analysis may be also applied in many other aspects.

Classical convex analysis is based on the framework of a locally convex space, so establishing random convex analysis requires a proper random generalization of the notion of a locally convex space.  It is well known that a locally convex space may be defined in the two equivalent ways: (1) a locally convex space is an ordered pair $(E,\mathcal{P})$ such that $E$ is a linear space and $\mathcal{P}$ a family of seminorms on $E$; (2) a locally convex space is an ordered pair $(E,\mathcal{T})$ such that $(E,\mathcal{T})$ is a linear topological space and $\mathcal{T}$ has a local base of convex neighborhoods of the zero element. Guo gave the notion of a random locally convex module in \cite{TXG-basic,TXG-Sur}, which is a proper random generalization of the above first definition. Whereas, as a proper random generalization of the above second definition, D. Filipovi\'{c}, M. Kupper and N. Vogelpoth gave the notion of a locally $L^{0}$-convex module in \cite{FKV}, in company of which the important notion of a locally $L^0$-convex topology is introduced. The two kinds of random generalizations are no longer equivalent and their relations are rather complicated, which makes the study of random convex analysis more involved than that of classical convex analysis.

Random convex analysis was first studied under the framework of a locally $L^{0}$--convex module in \cite{FKV}. Obviously, the work in \cite{FKV} heavily depends on the premise that the locally $L^{0}$--convex topology for every locally $L^{0}$--convex module can be induced by a family of $L^{0}$--seminorms, namely Theorem 2.4 of \cite{FKV}. Unfortunately, Wu and Guo \cite{Wu-Guo} and J.M.Zapata\cite{Zapata} recently, independently gave a counterexample showing Theorem 2.4 of \cite{FKV} is false, so a locally $L^{0}$--convex module is not a proper space framework for random convex analysis. In this paper we present a new approach to random convex analysis in order to overcome this obstacle by choosing a random locally convex module as the framework on which random convex analysis is based. Since the notion of a random locally convex module itself assumes the existence of a family of $L^{0}$--seminorms, furthermore, the notion of a random locally convex module has the advantage that it may be simultaneously endowed with the two kinds of topologies (namely the $(\varepsilon,\lambda)$--topology and the locally $L^0$-- convex topology) and the relations between some basic results derived from the two kinds of topologies have been given in \cite{TXG-JFA}, we can now present a complete random convex analysis.

According to our plan in \cite{GZZ1}, our work on random convex analysis can be divided into the three parts: (1) separation and Fenchel-Moreau duality in random locally convex modules; (2) continuity and subdifferentiability theorems in $L^{0}$--pre--barreled random locally convex modules; (3) the relations among the three kinds of $L^0$-convex conditional risk measures, cf. \cite{GZZ}. In order for this paper to be not too long, this paper is only focused on the above first part, in particular we establish the precise relation between random conjugate spaces of a random locally convex module under the two kinds of topologies and find the nice property of a local function on a countable concatenation hull of a set, which plays a crucial role in the establishing process of random convex analysis as well as in this paper. Besides, this paper also provides an important example (Example 4.11) exhibiting some pitfalls in the study of separation between a point and a closed $L^{0}$-convex subset in a random locally convex module. In fact, the results of this paper have been used in \cite{GZZ}.

The remainder of this paper is organized as follows. Section 2 recalls some necessary basic concepts; Section 3 establishes the precise relations between the two kinds of random conjugate spaces of a random locally convex module under the $(\varepsilon,\lambda)$--topology and locally $L^{0}$--convex topology; Section 4 gives a thorough treatment of separation between a point and a closed $L^{0}$-convex subset in a random locally convex module; Section 5 proves the complete Fenchel-Moreau duality theorems in random locally convex modules under the two kinds of topologies.

Throughout this paper, we always use the following notation and terminology:

$K:$ the scalar field R of real numbers or C of complex numbers.

$(\Omega,\mathcal{F},P):$ a probability space.

$L^{0}(\mathcal{F},K)=$ the algebra of equivalence classes of $K$--valued $\mathcal{F}$-- measurable random variables
on $(\Omega,\mathcal{F},P)$.

$L^{0}(\mathcal{F})=L^{0}(\mathcal{F},R)$.

$\bar{L}^{0}(\mathcal{F})=$ the set of equivalence classes of extended real-valued $\mathcal{F}$-- measurable random
variables on $(\Omega,\mathcal{F},P)$.

As usual, $\bar{L}^{0}(\mathcal{F})$ is partially ordered by $\xi\leq\eta$ iff $\xi^{0}(\omega)\leq\eta^{0}(\omega)$ for $P$--almost all $\omega\in \Omega$ (briefly, a.s.), where $\xi^0$ and $\eta^0$ are arbitrarily chosen representatives of $\xi$ and $\eta$, respectively. Then $(\bar{L}^{0}(\mathcal{F}),\leq)$ is a complete lattice, $\bigvee H$ and $\bigwedge H$ denote the supremum and infimum of a subset $H$, respectively. $(L^{0}(\mathcal{F}),\leq)$ is a conditionally complete lattice. Please refer to \cite{D-Sch} or \cite[p. 3026]{TXG-JFA} for the rich properties of the supremum and infimum of a set in $\bar{L}^{0}(\mathcal{F})$.

Let $\xi$ and $\eta$ be in $\bar{L}^{0}(\mathcal{F})$. $\xi<\eta$ is understood as usual, namely $\xi\leq\eta$ and $\xi\neq\eta$. In this paper we also use $``\xi<\eta $ (or $\xi\leq\eta$) on $A"$ for $``\xi^{0}(\omega)<\eta^{0}(\omega)$ (resp., $\xi^{0}(\omega)\leq\eta^{0}(\omega)$) for $P$--almost all $\omega\in A"$, where $A\in\mathcal{F}$, $\xi^0$ and $\eta^0$ are a representative of $\xi$ and $\eta$, respectively.

$\bar{L}^{0}_{+}(\mathcal{F})=\{\xi\in \bar{L}^{0}(\mathcal{F})~|~\xi\geq0\}$

$L^{0}_{+}(\mathcal{F})=\{\xi\in L^{0}(\mathcal{F})~|~\xi\geq0\}$

$\bar{L}^{0}_{++}(\mathcal{F})=\{\xi\in \bar{L}^{0}(\mathcal{F})~|~\xi>0$ on $\Omega\}$

$L^{0}_{++}(\mathcal{F})=\{\xi\in L^{0}(\mathcal{F})~|~\xi>0$ on $\Omega\}$

Besides, $\tilde{I}_{A}$ always denotes the equivalence class of $I_{A}$, where $A\in \mathcal{F}$ and $I_{A}$ is the characteristic function of $A$.
When $\tilde{A}$ denotes the equivalence class of $A (\in \mathcal{F})$, namely $\tilde{A}=\{B\in\mathcal{F}~|~P(A\triangle B)=0\}$ (here, $A\triangle B=(A\setminus B)\bigcup(B\setminus A)$), we also use $I_{\tilde{A}}$ for $\tilde{I}_{A}$.

Specially, $[\xi<\eta]$ denotes the equivalence class of $\{\omega\in\Omega~|~\xi^0(\omega)<\eta^0(\omega)\}$, where $\xi^0$ and $\eta^0$ are arbitrarily chosen representatives of $\xi$ and $\eta$ in $\bar{L}^0(\mathcal{F})$, respectively, some more notations such as $[\xi=\eta]$ and $[\xi\neq\eta]$ can be similarly understood.

\section{Some basic concepts}

Let us first recall the notion of a random normed module since its study motivates the notions of a random locally convex module and a locally $L^{0}$--convex module.

\begin{definition}($See$ \cite{TXG-Master,TXG-PHD,TXG-basic}). An ordered pair $(E,\|\cdot\|)$ is called a random normed space (briefly, an $RN$ space) over $K$ with base $(\Omega,\mathcal{F},P)$ if $E$ is a linear space over $K$ and $\|\cdot\|$ is a mapping from $E$ to $L^{0}_{+}(\mathcal{F})$ such that the following are satisfied:

\noindent ($RN$--1). $\|\alpha x\|=|\alpha| \|x\|$, $\forall \alpha\in K$ and $x\in E$;

\noindent ($RN$--2). $\|x\|=0$ implies $x=\theta$ (the null element of $E$);

\noindent ($RN$--3). $\|x+y\|\leq\|x\|+\|y\|$, $\forall x,y\in E$.

\noindent Here $\|\cdot\|$ is called the random norm on $E$ and $\|x\|$ the random norm of $x\in E$ (If $\|\cdot\|$ only satisfies ($RN$--1) and ($RN$--3) above, it is called a random seminorm on $E$).

\noindent Furthermore, if, in addition, $E$ is a left module over the algebra $L^{0}(\mathcal{F},K)$ (briefly, an $L^{0}(\mathcal{F},K)$--module) such that

\noindent ($RNM$--1). $\|\xi x\|=|\xi| \|x\|$, $\forall \xi\in L^{0}(\mathcal{F},K)$ and $x\in E$.

\noindent Then $(E,\|\cdot\|)$ is called a random normed module (briefly, an $RN$ module) over $K$ with base $(\Omega,\mathcal{F},P)$, the random norm $\|\cdot\|$ with the property ($RNM$--1) is also called an $L^0$--norm on $E$ (a mapping only satisfying ($RN$--3) and ($RNM$--1) above is called an $L^0$--seminorm on E).

\end{definition}

\begin{remark} According to the original notion of an $RN$ space in \cite{SS}, $\|x\|$ is a nonnegative random variable for all $x\in E$. An $RN$ space in the sense of Definition 2.1 is almost equivalent to (in fact, slightly more general than) the original one in the sense of \cite{SS}. Definition 2.1 is not only very natural from traditional functional analysis but also easily avoids any possible ambiguities between random variables and their equivalence classes, and hence also more convenient for applications to Lebesgue-Bochner function spaces since the latter exactly consists of equivalence classes. $RN$ spaces in the sense of Definition 2.1 was essentially earlier employed in \cite{TXG-Master}. The study of random conjugate spaces (see Definition 2.3 below) of $RN$ spaces and applications of $RN$ spaces to best approximations in Lebesgue-Bochner function spaces lead Guo to the notion of an $RN$ module in \cite{TXG-PHD}. Subsequently, $RN$ modules and their random conjugate spaces were deeply developed by Guo in \cite{TXG-PHD,TXG-Extension,TXG-Module,TXG-Radon} so that Guo further presented the refined notions of $RN$ modules and compared the original notion of $RN$ spaces with the currently used one of $RN$ spaces in \cite{TXG-basic}. At almost the same time, as a tool for the study of ultrapowers of Lebesgue-Bochner function spaces, $RN$ spaces and $RN$ modules were independently introduced by R. Haydon, M. Levy and Y. Raynaud in \cite{HLR}, where their notion of randomly normed $L^0(\mathcal{F})$--modules is exactly that of $RN$ modules over $R$ with base $(\Omega,\mathcal{F},P)$, in particular, they deeply studied the two classes of $RN$ modules-direct integrals and random Banach lattices (namely, random normed module equivalent of Banach lattices). Motivated by financial applications, D. Filipovi\'{c}, M. Kupper and N. Vogelpoth also independently came to the notion of $RN$ modules in \cite{FKV}, where the notion of $L^0$--normed modules amounts to that of $RN$ modules over $R$ with base $(\Omega,\mathcal{F},P)$.
\end{remark}

Following are the two important examples of $RN$ modules.

D. Filipovi$\acute{c}$, M. Kupper and N. Vogelpoth constructed important $RN$ modules $L^p_{\mathcal{F}}(\mathcal{E}) (1\leq p\leq+\infty)$ in \cite{FKV}, we will prove that they play the role of universal model spaces for $L^0$--convex conditional risk measures in the forthcoming paper.

\begin{example} Let $(\Omega, {\mathcal E}, P)$ be a probability space and ${\mathcal F}$
a sub--$\sigma$--algebra of ${\mathcal E}$. Define $|||\cdot|||_p\colon L^0({\mathcal E})\to {\bar L}^0_+({\mathcal F})$ by
$$|||x|||_p=\left\{
               \begin{array}{ll}
                 E[|x|^p|{\mathcal F}]^{1\over p}, & \hbox{when $1\leq p<\infty$;} \\
                 \bigwedge\{\xi\in {\bar L}^0_+({\mathcal F})~|~|x|\leq\xi\}, & \hbox{when $p=+\infty$;}
               \end{array}
             \right.
$$
for all $x\in L^0({\mathcal E})$.

Denote $L^p_{\mathcal F}({\mathcal E})=\{x\in L^0({\mathcal E})~|~|||x|||_p\in L^0_+({\mathcal F})\}$, then $(L^p_{\mathcal F}({\mathcal E}), |||\cdot|||_p)$ is an $RN$ module over $R$ with base $(\Omega, {\mathcal F}, P)$ and $L^p_{\mathcal F}({\mathcal E})=L^0({\mathcal F})\cdot L^p({\mathcal E})=\{~\xi x~|~\xi\in L^0({\mathcal F})~ \hbox{and}~ x\in L^p({\mathcal E})\}$.

\end{example}

To put some important classes of stochastic processes into the framework of $RN$ modules, Guo constructed a more general $RN$ module
$L^p_{\mathcal F}(S)$ in \cite{TXG-JFA} for each $p\in [1, +\infty]$, one can imagine that $S$ is an $RN$ module generated by a class of stochastic processes, $L^p_{\mathcal F}(S)$ can be constructed as follows.

\begin{example} Let $(S, \|\cdot\|)$ be an $RN$ module over $K$ with base $(\Omega, {\mathcal E}, P)$ and ${\mathcal F}$
a sub--$\sigma$--algebra. Define $|||\cdot|||_p\colon S \to {\bar L}^0_+({\mathcal F})$ by
$$|||x|||_p=\left\{
               \begin{array}{ll}
                 E[\|x\|^p|{\mathcal F}]^{1\over p}, & \hbox{when $1\leq p<\infty$;} \\
                 \bigwedge\{\xi\in {\bar L}^0_+({\mathcal F})|~\|x\|\leq\xi\}, & \hbox{when $p=+\infty$;}
               \end{array}
             \right.
$$
for all $x\in S$.

Denote $L^p_{\mathcal F}(S)=\{x\in S~|~|||x|||_p\in L^0_+({\mathcal F})\}$, then $(L^p_{\mathcal F}(S), |||\cdot|||_p)$ is an $RN$ module over $K$ with base $(\Omega, {\mathcal F}, P)$. When $S=L^0({\mathcal E})$, $L^p_{\mathcal F}(S)$ is exactly $L^p_{\mathcal F}({\mathcal E})$.
\end{example}

\begin{definition}($See$ \cite{TXG-PHD,TXG-Module,TXG-Sur}). An ordered pair $(E,\mathcal{P})$ is called a random locally convex space (briefly, an $RLC$ space) over $K$ with base $(\Omega,\mathcal{F},P)$ if $E$ is a linear space over $K$ and $\mathcal{P}$ a family of mappings from $E$ to $L^0_{+}(\mathcal{F})$ such that the following are satisfied:

\noindent ($RLC$--1). Every $\|\cdot\|\in \mathcal{P}$ is a random seminorm on $E$;

\noindent ($RLC$--2). $\bigvee\{\|x\|:\|\cdot\|\in\mathcal{P}\}=0$ iff $x=\theta$.

\noindent Furthermore, if, in addition, $E$ is an $L^{0}(\mathcal{F},K)$--module and each $\|\cdot\|\in \mathcal{P}$ is an $L^0$--seminorm on $E$, then $(E,\mathcal{P})$ is called a random locally convex module (briefly, an $RLC$ module) over $K$ with base $(\Omega,\mathcal{F},P)$.
\end{definition}

In the sequel of this paper, given a random locally convex space $(E,\mathcal{P})$, $\mathcal{P}_f$ always denotes the family of finite subsets of $\mathcal{P}$, for each $\mathcal{Q}\in \mathcal{P}_f$ $\|\cdot\|_{\mathcal{Q}}$ denotes the random seminorm defined by $\|x\|_{\mathcal{Q}}=\bigvee\{\|x\|:\|\cdot\|\in \mathcal{Q}\}$ for all $x\in E$ and $\mathcal{P}_{cc}=\{\sum_{n=1}^{\infty}\tilde{I}_{A_n}\|\cdot\|_{\mathcal{Q}_n}~|~\{A_n,n\in N\}$ is a countable partition of $\Omega$ to $\mathcal{F}$ and $\{\mathcal{Q}_n,n\in N\}$ a sequence of finite subsets of $\mathcal{P}\}$, called the countable concatenation hull of $\mathcal{P}$.

\begin{definition}($See$ \cite{TXG-PHD,TXG-Module,TXG-Sur,TXG-SLP}). Let $(E, {\mathcal P})$ be an $RLC$ space over $K$ with base $(\Omega, {\mathcal F}, P)$. For any positive numbers $\varepsilon$ and $\lambda$ with $0<\lambda<1$ and $\mathcal{Q}\in {\mathcal P}_f$, let $N_{\theta}(\mathcal{Q}, \varepsilon, \lambda)=\{x\in E~|~P\{\omega\in \Omega~|~\|x\|_\mathcal{Q}(\omega)<\varepsilon\}>1-\lambda\}$, then $\{N_{\theta}(\mathcal{Q}, \varepsilon, \lambda)~|~\mathcal{Q}\in {\mathcal P}_f, \varepsilon >0, 0<\lambda<1\}$ forms a local base at $\theta$ of some Hausdorff linear topology on $E$, called the $(\varepsilon, \lambda)$--topology induced by ${\mathcal P}$.
\end{definition}

From now on, we always denote by ${\mathcal T}_{\varepsilon, \lambda}$ the $(\varepsilon, \lambda)$--topology for every $RLC$ space if there is no possible confusion. Clearly, the $(\varepsilon, \lambda)$--topology for the special class of $RN$ modules $L^0({\mathcal F}, K)$ is exactly the ordinary topology of convergence in measure, and $(L^0({\mathcal F}, K), {\mathcal T}_{\varepsilon, \lambda})$ is a topological algebra over $K$. It is also easy to check that $(E, {\mathcal T}_{\varepsilon, \lambda})$ is a topological module over $(L^0({\mathcal F}, K), {\mathcal T}_{\varepsilon, \lambda})$ when $(E, {\mathcal P})$ is an $RLC$ module over $K$ with base $(\Omega, {\mathcal F}, P)$, namely the module multiplication operation is jointly continuous.

For any $\varepsilon \in L^0_{++}({\mathcal F})$, let $U(\varepsilon)=\{\xi\in L^0({\mathcal F}, K)~|~|\xi|\leq \varepsilon\}$. A subset $G$ of $L^0({\mathcal F}, K)$ is ${\mathcal T}_c$--open if for each fixed $x\in G$ there is some $\varepsilon \in L^0_{++}({\mathcal F})$ such that $x+U(\varepsilon)\subset G$. Denote by ${\mathcal T}_c$ the family of ${\mathcal T}_c$--open subsets of $L^0({\mathcal F}, K)$, then ${\mathcal T}_c$ is a Hausdorff topology on $L^0({\mathcal F}, K)$ such that $(L^0({\mathcal F}, K), {\mathcal T}_c)$ is a topological ring, namely the addition and multiplication operations are jointly continuous. D. Filipovi\'{c}, M. Kupper and N. Vogelpoth first observed this kind of topology and further pointed out that ${\mathcal T}_c$ is not necessarily a linear topology since the mapping $\alpha\mapsto \alpha x$ ($x$ is fixed) is no longer continuous in general. These observations led them to the study of a class of topological modules over the topological ring $(L^0({\mathcal F}, K), {\mathcal T}_c)$ in \cite{FKV}, where they only considered the case when $K=R$, in fact the complex case can also similarly introduced as follows.

\begin{definition}($See$ \cite{FKV}). An ordered pair $(E, {\mathcal T})$ is a topological $L^0({\mathcal F}, K)$--module if both $(E, {\mathcal T})$ is a topological space and $E$ is an $L^0({\mathcal F}, K)$--module such that $(E, {\mathcal T})$ is a topological module over the topological ring $(L^0({\mathcal F}, K), {\mathcal T}_c)$, namely the addition and module multiplication operations are jointly continuous.
\end{definition}

\begin{definition}($See$ \cite{TXG-Sur,TXG-XXC,FKV}). Let $E$ be an $L^0({\mathcal F}, K)$--module and $A$ and $B$ two subsets of $E$. $A$ is said to be $L^0$--absorbed by $B$ if there is some $\xi \in L^0_{++}({\mathcal F})$ such that $\eta A\subset B$ for all $\eta\in L^0({\mathcal F}, K)$ with $|\eta|\leq \xi$. $B$ is $L^0$--absorbent if $B$ $L^0$--absorbs every element in $E$. $B$ is $L^0$--convex if $\xi x+(1-\xi)y\in B$ for all $x,\,y\in B$ and $\xi\in L^0_+({\mathcal F})$ with $0\leq \xi\leq 1$. $B$ is $L^0$--balanced if $\eta B\subset B$ for all $\eta\in L^0({\mathcal F}, K)$ with $|\eta|\leq 1$.
\end{definition}

\begin{remark} Clearly, when $B$ is $L^0$--balanced, $A$ is $L^0$--absorbed by $B$ iff there exists some $\xi \in L^0_{++}({\mathcal F})$ such that $A\subset \xi B$. Since $L^0({\mathcal F}, K)$ is an algebra over $K$, an $L^0({\mathcal F}, K)$--module is also a linear space over $K$, then it is clear that $B$ is balanced (resp., convex) if $B$ is $L^0$--balanced (resp., $L^0$--convex). But `` being $L^0$--absorbent '' and `` being absorbent '' may not imply each other.
\end{remark}

\begin{definition}($See$ \cite{FKV}). A topological $L^0({\mathcal F}, K)$--module $(E, {\mathcal T})$ is called a locally $L^0$--convex $L^0({\mathcal F}, K)$--module ( briefly, a locally $L^0$--convex module when $K=R$ ), in which case ${\mathcal T}$ is called a locally $L^0$--convex topology on $E$, if ${\mathcal T}$ has a local base ${\mathcal B}$ at $\theta$ ( the null element in $E$ ) such that each member in ${\mathcal B}$ is $L^0$--balanced, $L^0$--absorbent and $L^0$--convex.
\end{definition}

\begin{proposition}($See$ \cite{FKV}). Let ${\mathcal P}$ be a family of $L^0$--seminorms on an $L^0({\mathcal F}, K)$--module $E$. For any $\varepsilon \in L^0_{++}({\mathcal F})$ and any $ Q\in {\mathcal P}_f$ (namely $Q$  is a finite subset of ${\mathcal P}$), let $N_{\theta}(Q, \varepsilon)=\{x\in E~|~\|x\|_Q\leq \varepsilon\}$, then $\{~N_{\theta}(Q, \varepsilon)~|~Q\in {\mathcal P}_f, ~\varepsilon \in L^0_{++}({\mathcal F})\}$ forms a local base at $\theta$ of some locally $L^0$--convex topology, called the locally $L^0$--convex topology induced by ${\mathcal P}$.
\end{proposition}

\begin{corollary} Let $(E,{\mathcal P})$ be an $RLC$ module over $K$ with base $(\Omega, {\mathcal F}, P)$ and ${\mathcal T}_c$ the locally $L^0$--convex topology induced by ${\mathcal P}$. Then $(E, {\mathcal T}_c)$ is a Hausdorff locally $L^0$--convex $L^0({\mathcal F}, K)$--module.
\end{corollary}

From now on, we always denote by $\mathcal{T}_c$ the locally $L^0$--convex topology induced by $\mathcal{P}$ for every $RLC$ module $(E,\mathcal{P})$ if there is no risk of confusion.

In the final part of this section, let us return to the basic problem: whether can a locally $L^0$--convex topology on an $L^0({\mathcal F}, K)$--module $E$ be induced by a family of $L^0$--seminorms on $E$? If the answer is yes, then the theory of Hausdorff locally $L^0$--convex modules is equivalent to that of random locally convex modules endowed with the locally $L^0$--convex topology, which will be a perfect counterpart of the classical result that a Hausdorff locally convex topology can be induced by a separating family of seminorms. It is well known that classical gauge functionals play a crucial role in the proof of the classical result. Let $U$ be a balanced, absorbent and convex subset of a locally convex space $(E, {\mathcal T})$ and $p_U$ the gauge functional of $U$, then the following relation is easily verified: $$\{x\in E~|~p_U(x)<1\}\subset U\subset \{x\in E~|~p_U(x)\leq 1\}, \eqno(2.1)$$
It is the relation (2.1) that is key in the proof of the above classical result.

Random gauge functional was first introduced in \cite{FKV}. Let $U$ be an $L^0$--balanced, $L^0$--absorbent and $L^0$--convex subset of an $L^0({\mathcal F}, K)$--module $E$, define $p_U\colon E\to L^0_+({\mathcal F})$ by $p_U(x)=\bigwedge\{\xi\in L^0_+({\mathcal F})~|~x\in \xi U\}$ for all $x\in E$, called the random gauge functional of $U$. Furthermore, it is also proved in \cite{FKV} that $p_U(x)=\bigwedge\{\xi\in L^0_{++}({\mathcal F})~|~x\in \xi U\}$ for all $x\in E$ and $p_U$ is an $L^0$--seminorm on $E$.

Let $(E, {\mathcal T})$ be a locally $L^0$--convex $L^0({\mathcal F}, K)$--module and $\mathcal{U}$ a local base at the null of ${\mathcal T}$ such that each $U\in \mathcal{U}$ is $L^0$--balanced, $L^0$--absorbent and $L^0$--convex. Furthermore, let $\mathcal{P}=\{p_U~:~U\in \mathcal{U}\}$. If it was proved that $ \{x\in E~|~p_U(x)< 1~\hbox{on}~\Omega\}\subset U\subset\{x\in E~|~p_U(x)\leq1~\hbox{on}~\Omega\}$ for each $U\in \mathcal{U}$, then ${\mathcal T}$ would be equivalent to the locally $L^0$--convex topology induced by $\mathcal{P}$, namely every locally $L^0$--convex topology could be induced by a family of $L^0$--seminorms, that is to say that Theorem 2.4 of \cite{FKV} would be true. Unfortunately, Wu and Guo \cite{Wu-Guo} and J.M.Zapata\cite{Zapata} recently, independently constructed a counterexample showing that Theorem 2.4 of \cite{FKV} is false. In fact, D. Filipovi\'{c}, M. Kupper and N. Vogelpoth only proved that the locally $L^0$--convex topology induced by $\mathcal{P}$ is weaker than ${\mathcal T}$ in \cite{FKV} since they only proved the following:

\begin{proposition}($See$ \cite{FKV}). Let $(E, {\mathcal T})$ be a locally $L^0$--convex $L^0({\mathcal F}, K)$--module and $U$ an $L^0$--balanced, $L^0$--absorbent and $L^0$--convex subset of $E$. Then the following statements hold:\\
(i).~~$U\subset \{x\in E~|~p_U(x)\leq 1\}$;\\
(ii).~~$p_U(x)\geq 1 $ on $B$ if ${\tilde I}_Ax\notin {\tilde I}_AK$ for all $A\in {\mathcal F}$ with $P(A)>0$ and $A\subset B$, where $B\in {\mathcal F}$ satisfies $P(B)>0$;\\
(iii).~~$U^o\subset \{x\in E~|~p_U(x)<1~\hbox{on}~\Omega\}$.
\end{proposition}

The most interesting part in Proposition 2.13 is (ii). In fact, (i) is clear and (iii) can be proved as follows: Given an $x\in U^o$, there is an $L^0$--balanced, $L^0$--absorbent and $L^0$--convex neighborhood $V$ of $\theta$ such that $x+V\in U$. Since there is $\delta \in L^0_{++}({\mathcal F})$ such that $\delta x\in V$, $(1+\delta)x=x+\delta x\in x+V\subset U$, so $x\in {\frac{1}{1+\delta}}U$, then $p_U(x)\leq {\frac{1}{1+\delta}}<1~\hbox{on}~\Omega$.

Then, can (ii) of Proposition 2.13 imply that $\{x\in E~|~p_U(x)< 1~\hbox{on}~\Omega\}\subset U$? Or, we can ask: does it hold that $\{x\in E~|~p_U(x)< 1~\hbox{on}~\Omega\}\subset U$? Proposition 2.15 below shows that it is not a simple problem whether $\{x\in E~|~p_U(x)< 1~\hbox{on}~\Omega\}$ is contained in $U$.

Let us first recall the notion of countable concatenation property of a set or an $L^0({\mathcal F}, K)$--module. The introducing of the notion utterly results from the study of the locally $L^0$--convex topology, the reader will see that this notion is ubiquitous in the theory of the locally $L^0$--convex topology.

From now on, we always suppose that all the $L^0({\mathcal F}, K)$--modules $E$ involved in this paper have the property that for any $x,~y\in E$, if there is a countable partition $\{A_n,n\in N\}$ of $\Omega$ to ${\mathcal F}$ such that ${\tilde I}_{A_n}x={\tilde I}_{A_n}y$ for each $n\in N$ then $x=y$. Guo already pointed out in \cite{TXG-JFA} that all random locally convex modules possess this property, so the assumption is not too restrictive.

\begin{definition}($See$ \cite{TXG-JFA}). Let $E$ be an $L^0({\mathcal F}, K)$--module. A sequence $\{x_n, n\in N\}$ in $E$ is countably concatenated in $E$ with respect to a countable partition $\{A_n,n\in N\}$ of $\Omega$ to ${\mathcal F}$ if there is $x\in E$ such that ${\tilde I}_{A_n}x={\tilde I}_{A_n}x_n$ for each $n\in N$, in which case we define $\sum^{\infty}_{n=1}{\tilde I}_{A_n}x_n$ as $x$. A subset $G$ of $E$ is said to have the countable concatenation property if each sequence $\{x_n, n\in N\}$ in $G$ is countably concatenated in $E$ with respect to an arbitrary countable partition $\{A_n,n\in N\}$ of $\Omega$ to $\mathcal{F}$ and $\sum^{\infty}_{n=1}{\tilde I}_{A_n}x_n\in G$.
\end{definition}

From now on, let $E$ be an $L^0({\mathcal F}, K)$--module with the countable concatenation property and $G$ a subset of $E$. $H_{cc}(G)$ always denotes the countable concatenation hull of $G$ in $E$, namely $H_{cc}(G)=\{\sum^{\infty}_{n=1}{\tilde I}_{A_n}g_n:$ $\{A_n,n\in N\}$ is a countable partition of $\Omega$ to $\mathcal{F}$ and  $\{g_n,n\in N\}$ is a sequence in $G\}$. Furthermore, if $x=\sum^{\infty}_{n=1}{\tilde I}_{A_n}x_n$ for some countable partition $\{A_n,n\in N\}$ of $\Omega$ to ${\mathcal F}$ and some sequence $\{x_n,n\in N\}$ in $E$, then $\sum^{\infty}_{n=1}{\tilde I}_{A_n}x_n$ is called a canonical representation of $x$.

\begin{proposition} Let $(E, {\mathcal T})$ be a locally $L^0$--convex $L^0({\mathcal F}, K)$--module and $U$ an $L^0$--balanced, $L^0$--absorbent and $L^0$--convex subset with the countable concatenation property. Then $U^o\subset \{x\in E~|~p_U(x)<1 ~\hbox{on}~\Omega\}\subset U\subset \{x\in E~|~p_U(x)\leq 1\}$, where $U^o$ denotes the ${\mathcal T}$--interior of $U$.
\end{proposition}

\begin{proof} By Proposition 2.13, we only need to show that $\{x\in E~|~p_U(x)<1 ~\hbox{on}~\Omega\}\subset U$. Let $x_{0}$ be a point in $E$ such that $p_U(x_{0})<1~\hbox{on}~\Omega$. Since $\{\xi\in L^0_{++}({\mathcal F})~|~x_{0}\in \xi U\}$ is downward directed, there is a sequence $\{\xi_n, n\in N\}$ in $L^0_{++}({\mathcal F})$ such that it converges to $p_U(x_{0})$ in a nonincreasing way and $x_{0}\in \xi_n U$ for each $n\in N$. By the Egoroff theorem there are a countable partition $\{A_n,n\in N\}$ of $\Omega$ to ${\mathcal F}$ and a subsequence $\{\xi_{n_k},k\in N\}$ of $\{\xi_n, n\in N\}$ such that the subsequence converges to $p_U(x_{0})$ uniformly on each $A_n$. So, we can suppose that the subsequence is just $\{\xi_n, n\in N\}$ itself and each $\xi_n<1~\hbox{on}~A_n$ since $p_U(x_{0})<1 ~\hbox{on}~\Omega$. Clearly, ${\tilde I}_{A_n}x_{0}\in {\tilde I}_{A_n}\xi_n U$ for each $n\in N$, let $u_n\in U$ be such that ${\tilde I}_{A_n}x_{0}={\tilde I}_{A_n}\xi_n u_n$ for each $n\in N$. Let $u^\prime_n={\tilde I}_{A_n}\xi_n u_n$, then $u^\prime_n\in U$, it is obvious that the sequence $\{u^\prime_n, n\in N\}$ is countably concatenated with respect to $\{A_n, n\in N\}$, so that $x_0\in U$ since $U$ has the countable concatenation property.
\end{proof}

Proposition 2.15 tells us that for a locally $L^0$--convex $L^0({\mathcal F}, K)$--module $(E, {\mathcal T})$, if ${\mathcal T}$ has a local base consisting of $L^0$--balanced, $L^0$--absorbent and $L^0$--convex subsets with the countable concatenation property, then ${\mathcal T}$ can be induced by a family of $L^0$--seminorms, however, it is rather restrictive to require the existence of such a local base.
 It is Proposition 2.15 that motivates Wu and Guo \cite{Wu-Guo} and J.M.Zapata \cite{Zapata} to go farther. Precisely speaking, Wu and Guo \cite{Wu-Guo} and J.M.Zapata \cite{Zapata}, independently, introduced the notion of the relative countable concatenation property, which is weaker than that of the countable concatenation property but meets the needs of Wu and Guo \cite{Wu-Guo} and J.M.Zapata \cite{Zapata}, in fact, Wu and Guo \cite{Wu-Guo} and J.M.Zapata \cite{Zapata}, independently, have given a necessary and sufficient condition for a locally $L^0$--convex topology to be induced by a family of $L^0$--seminorms and in particular have given a counterexample showing that not every locally $L^0$--convex topology is necessarily induced by a family of $L^0$--seminorms.

\section{The precise relation between the random conjugate spaces of a random locally convex module under the two kinds of topologies}

The main result of this section is Theorem 3.7. Let us first recall the notion of the random conjugate space of an $RN$ space.

\begin{definition}($See$ \cite{TXG-Master,TXG-PHD,TXG-Module,TXG-basic}). Let $(E,\|\cdot\|)$ be an $RN$ space over $K$ with base $(\Omega,\mathcal{F},P)$. A linear operator $f$ from $E$ to $L^{0}(\mathcal{F},K)$ is said to be an a.s. bounded random linear functional if there is $\xi\in L^{0}_{+}(\mathcal{F})$ such that $\|f(x)\|\leq\xi\|x\|, \forall x\in E$. Denote by $E^{\ast}$ the linear space of a.s. bounded random linear functionals on $E$, define $\|\cdot\|:E^{\ast}\rightarrow L^{0}_{+}(\mathcal{F})$ by $\|f\|=\bigwedge\{\xi\in L^{0}_{+}(\mathcal{F})~|~\|f(x)\|\leq\xi\|x\|$ for all $x\in E\}$ for all $f\in E^{\ast}$, then it is easy to check that $(E^{\ast},\|\cdot\|)$ is also an $RN$ module over $K$ with base $(\Omega,\mathcal{F},P)$, called the random conjugate space of $E$.
\end{definition}

It is easy to see that Definition 3.1 coincides with the notion of a random dual introduced in \cite{HLR} of an $RN$ space. It is not very difficult to introduce the random conjugate space for an $RN$ space, whereas it is completely another thing to do for an $RLC$ space, at the earlier time Guo ever gave two definitions, which turns out to be equivalent to the two kinds of random conjugate spaces for a random locally convex module under the two kinds of topologies, see Propositions 3.3 and 3.4 below.

\begin{definition}($See$ \cite{TXG-PHD,TXG-Module,TXG-Sur}). Let $(E,\mathcal{P})$ be an $RLC$ space over $K$ with base $(\Omega,\mathcal{F},P)$. A linear operator $f$ from $E$ to $L^0(\mathcal{F},K)$ (such an operator is also called a random linear functional on $E$) is called an a.s. bounded random linear functional of type I if there are $\xi\in L^0_+(\mathcal{F})$ and some finite subset $\mathcal{Q}$ of $\mathcal{P}$ such that $|f(x)|\leq\xi \|x\|_{\mathcal{Q}}$ for all $x\in E$. Denote by $E^{\ast}_{I}$ the $L^0(\mathcal{F},K)$--module of a.s. bounded random linear functionals on $E$ of type I, called the first kind of random conjugate space of $(E,\mathcal{P})$, cf. \cite{TXG-PHD,TXG-Module}. A random linear functional $f$ on $E$ is called an a.s. bounded random linear functional of type II if there are $\xi\in L^0_+(\mathcal{F})$ and $\|\cdot\|\in \mathcal{P}_{cc}$ such that $|f(x)|\leq\xi\|x\|$ for all $x\in E$. Denote by $E^{\ast}_{II}$ the $L^0(\mathcal{F},K)$--module of a.s. bounded random linear functionals on $E$ of type II, called the second kind of random conjugate space of $(E,\mathcal{P})$, cf. \cite{TXG-Sur}.
\end{definition}

For an $RLC$ module $(E, {\mathcal P})$ over $K$ with base $(\Omega, {\mathcal F}, P)$, we always denote by $(E, {\mathcal P})^\ast_{\varepsilon, \lambda}$ ( or, briefly, $E^\ast_{\varepsilon, \lambda}$, whenever there is no confusion ) the $L^0({\mathcal F}, K)$--module of continuous module homomorphisms from $(E, {\mathcal T}_{\varepsilon, \lambda})$ to $(L^0({\mathcal F}, K), {\mathcal T}_{\varepsilon, \lambda})$, called the random conjugate space of $(E, {\mathcal P})$ under the  $(\varepsilon, \lambda)$--topology.

Guo proved that a linear operator $f$ from an $RN$ module $(E,\|\cdot\|)$ over $K$ with base $(\Omega, {\mathcal F}, P)$ to $L^0({\mathcal F}, K)$ is a.s. bounded if and only if $f$ is a continuous module homomorphism from $(E,\mathcal{T}_{\varepsilon, \lambda})$ to $(L^0({\mathcal F}, K),\mathcal{T}_{\varepsilon, \lambda})$, namely $E^\ast=E^\ast_{\varepsilon, \lambda}$ for every $RN$ module $E$, cf. \cite{TXG-PHD,TXG-Extension}. This can be extended to the following more general case when $E$ is an $RLC$ module.

\begin{proposition}($See$ \cite{TXG-Sur,TXG-LHZ}). Let $(E, {\mathcal P})$ be an $RLC$ module $(E, {\mathcal P})$ over $K$ with base $(\Omega, {\mathcal F}, P)$ and $f$ a random linear functional on $E$. Then $f\in E^\ast_{II}$ iff $f\in E^\ast_{\varepsilon, \lambda}$, namely $E^\ast_{II}=E^\ast_{\varepsilon, \lambda}$.
\end{proposition}

Denote by $(E, {\mathcal T})^\ast_c$ ( briefly, $E^\ast_c$ ) the $L^0({\mathcal F}, K)$--module of continuous module homomorphisms from $(E, {\mathcal T})$ to $(L^0({\mathcal F}, K), {\mathcal T}_c)$, called the random conjugate space of the topological $L^0({\mathcal F}, K)$--module $(E, {\mathcal T})$, which was first introduced in \cite{FKV}.

Let $(E,{\mathcal P})^\ast_c=(E,{\mathcal T}_c)^\ast_c$ (briefly, $E^\ast_c$, if there is no risk of confusion), called the random conjugate space of a random locally convex module $(E,{\mathcal P})$ under the locally $L^0$--convex topology ${\mathcal T}_c$ induced by ${\mathcal P}$.

\begin{proposition}($See$ \cite{TXG-JFA}). Let $(E,{\mathcal P})$ be a random locally convex module over $K$ with base $(\Omega, {\mathcal F}, P)$ and $f\colon E\to L^0({\mathcal F}, K)$ a random linear functional. Then $f\in E^\ast_I$ iff $f\in E^\ast_c$, namely $E^\ast_I=E^\ast_c$.
\end{proposition}

\begin{remark} In \cite{TXG-JFA}, it is proved that $ E^\ast_c\subset E^\ast_I$ ( see \cite[p.3032]{TXG-JFA} ). Conversely, if $f\in E^\ast_I$, namely $f$ is a random linear functional and there are some $\xi\in L^0_+({\mathcal F})$ and $Q\in {\mathcal P}_f$ such that $|f(x)|\leq \xi \|x\|_Q$ for all $x\in E$. Lemma 2.12 of \cite{TXG-JFA} shows that $f$ must be $L^0({\mathcal F}, K)$--linear. It is also clear that $f$ is continuous from $(E,{\mathcal T}_c)$ to $(L^0({\mathcal F}, K), {\mathcal T}_c)$, and hence $f\in E^\ast_c$. Thus Proposition 3.4 was already proved in \cite{TXG-JFA}.
\end{remark}

A family ${\mathcal P}$ of random seminorms on a linear space $E$ is said to have the countable concatenation property if ${\mathcal P}_{cc}={\mathcal P}$, this definition appears stronger than that given in \cite{FKV} for a family of $L^0$--seminorms on an $L^0({\mathcal F}, K)$--module since ${\mathcal P}$ must be invariant under the operation of finitely many suprema once ${\mathcal P}_{cc}={\mathcal P}$. But ${\mathcal P}$ and $\{~\|\cdot\|_Q:Q\in {\mathcal P}_f\}$ always induces the same locally $L^0$--convex topology for any family ${\mathcal P}$ of $L^0$--seminorms on an $L^0({\mathcal F}, K)$--module $E$, thus the definition is essentially equivalent to that introduced in \cite{FKV}. It is also obvious that $E^\ast_I=E^\ast_{II}$ if ${\mathcal P}_{cc}={\mathcal P}$, and hence we have the following:

\begin{corollary} ($See$ \cite{TXG-JFA}). Let $(E,{\mathcal P})$ be an $RLC$ module over $K$ with base $(\Omega, {\mathcal F}, P)$. Then $E^\ast_c=E^\ast_{\varepsilon, \lambda}$ if ${\mathcal P}$ has the countable concatenation property. Specially, $E^\ast=E^\ast_c=E^\ast_{\varepsilon, \lambda}$ for an $RN$ module $(E,\|\cdot\|)$.
\end{corollary}

For an $RLC$ module $(E, {\mathcal P})$, by definition we have that $E^\ast_I\subset E^\ast_{II}$, so $E^\ast_c\subset E^\ast_{\varepsilon, \lambda}$ by Propositions 3.3 and 3.4. Guo pointed out in \cite{TXG-JFA} that $E^\ast_{\varepsilon, \lambda}$ has the countable concatenation property, if we denote by $H_{cc}(E^{\ast}_{c})$ the countable concatenation hull of $E^{\ast}_{c}$ in $E^\ast_{\varepsilon, \lambda}$, then we have the main result of this section, namely Theorem 3.7 below, which gives the precise relation between $E^\ast_c$ and $E^\ast_{\varepsilon, \lambda}$ for any $RLC$ module $E$.

\begin{theorem} Let $(E, {\mathcal P})$ be an $RLC$ module over $K$ with base $(\Omega, {\mathcal F}, P)$. Then  $E^{\ast}_{\varepsilon,\lambda}=H_{cc}(E^{\ast}_{c})$, where $E^{\ast}_{\varepsilon,\lambda}=(E, {\mathcal P})^\ast_{\varepsilon, \lambda}$ and $E^{\ast}_{c}=(E,\mathcal{P})^{\ast}_{c}$.
\end{theorem}

To prove Theorem 3.7, we first recall Lemma 3.8 below from \cite{TXG-XXC}.

\begin{lemma}($See$ \cite{TXG-XXC}). Let $X$ be an $L^{0}(\mathcal{F},K)$--module, $f:X\rightarrow L^{0}(\mathcal{F},K)$ an $L^{0}(\mathcal{F},K)$--linear function, $\{p_n:X\rightarrow L^{0}_{+}(\mathcal{F})~|~n\in N\}$ a sequence of $L^0$--seminorms on $X$ and $\{A_n,n\in N\}$ a countable partition of $\Omega$ to $\mathcal{F}$ such that $|f(x)|\leq\Sigma_{n=1}^{\infty}\tilde{I}_{A_n}p_n(x)$ for all $x\in X$. Then there is a sequence $\{f_n:n\in N\}$ of $L^{0}(\mathcal{F},K)$--linear functions such that

\noindent (1). $|f_n(x)|\leq p_n(x)$ for all $x\in X$  and $n\in N$;

\noindent (2). $f(x)=\Sigma_{n=1}^{\infty}\tilde{I}_{A_n}(f_n(x))$ for all $x\in X$.
\end{lemma}

We can now prove Theorem 3.7.

\noindent{\bf Proof of Theorem 3.7.}
Since $\mathcal{P}$ and $\mathcal{P}_{cc}$ induces the same $(\varepsilon,\lambda)$--topology on $E$, $E^{\ast}_{\varepsilon,\lambda}:=(E,\mathcal{P})^{\ast}_{\varepsilon,\lambda}=(E,\mathcal{P}_{cc})^{\ast}_{\varepsilon,\lambda}=(E,\mathcal{P}_{cc})^{\ast}_{c}$, where the last equality comes from the countable concatenation property of $\mathcal{P}_{cc}$ by Corollary 3.6. It remains to prove that $(E,\mathcal{P}_{cc})^{\ast}_{c}=H_{cc}(E^{\ast}_c)$ and we only needs to check that $(E,\mathcal{P}_{cc})^{\ast}_{c}\subset H_{cc}(E^{\ast}_c)$.

Let $f$ be any element of $(E,\mathcal{P}_{cc})^{\ast}_{c}$. Since $\mathcal{P}_{cc}$ is invariant under the operation of finitely many suprema, then there are $\|\cdot\|\in \mathcal{P}_{cc}$ and $\xi\in L^{0}_{++}(\mathcal{F})$ such that $|f(x)|\leq\xi\|x\|$ for all $x\in E$. Let $\|\cdot\|=\Sigma_{n=1}^{\infty}\tilde{I}_{A_n}\|\cdot\|_{\mathcal{Q}_n}$, where $\{A_n,n\in N\}$ is some countable partition of $\Omega$ to $\mathcal{F}$ and each $\mathcal{Q}_n\in \mathcal{P}_f$, then by Lemma 3.8 there is a sequence $\{f_n,n\in N\}$ of $L^{0}(\mathcal{F},K)$--linear functions such that

\noindent (1). $|f_n(x)|\leq\xi\|x\|_{\mathcal{Q}_n}$ for all $x\in E$ and $n\in N$;

\noindent (2). $f(x)=\Sigma_{n=1}^{\infty}\tilde{I}_{A_n}(f_n(x))$ for all $x\in E$.

(1) shows that each $f_n\in E^{\ast}_c$ and (2) further shows that $f=\Sigma_{n=1}^{\infty}\tilde{I}_{A_n}f_n$, so $f\in H_{cc}(E^{\ast}_c)$.
\hfill $\square$

\vspace{0.2cm}
For the further study of random conjugate spaces, please refer to \cite{TXG-SEZ}.

\section{Separation between a point and a closed $L^0$-convex subset in a random locally convex module}

As the classical hyperplane separation theorem between a point and a closed convex set in a locally convex space plays an essential role in the proof of the classical Fenchel-Moreau duality theorem, the hyperplane separation theorem between a point and a closed $L^0$-convex subset in a random locally convex module was also studied in \cite{TXG-strict,FKV,TXG-JFA}. Since there are the two kinds of topologies and correspondingly the two kinds of random conjugate spaces for a random locally convex module, the hyperplane separation theorem between a point and a closed $L^0$-convex subset in a random locally convex module also has the two forms, which corresponds to the two kinds of topologies, respectively. The aim of this section is to give some variants and improved versions of the corresponding results of \cite{TXG-JFA} and \cite{FKV} for the further study of random convex analysis as well as the proof of the main results of Section 5.

\subsection{Separation under the $(\varepsilon,\lambda)$-topology}

Let $(E, {\mathcal P})$ be an random locally convex module over $K$ with base $(\Omega, {\mathcal F}, P)$, $x\in E$ and $M$ a subset of $E$.  Define $d(x, M)=\bigvee\{d_{\mathcal{Q}}(x, M): \mathcal{Q}\in {\mathcal P}_f\}$, where $d_{\mathcal{Q}}(x, M)=\bigwedge\{\|x-y\|_{\mathcal{Q}}: y\in M\}$ for $\mathcal{Q}\in {\mathcal P}_f$. In this paper, we always use $(d(x, M)>0)$ for any chosen representative of $[d(x, M)>0]$. When $M$ is a $\mathcal{T}_{\varepsilon,\lambda}$-closed subset such that $\tilde{I}_{A}M+\tilde{I}_{A^c}M\subset M$, (3) of Lemma 3.8 of \cite{TXG-JFA} shows that $x\not\in M$ iff $d(x, M)>0$ (namely $(d(x, M)>0)$ has positive probability).

Proposition 4.1 below is due to \cite{TXG-strict}, the current form has been implied in the process of the proof of Theorem 3.7 of \cite{TXG-JFA}.

\begin{proposition} ($See$ \cite{TXG-JFA}). Let $(E, {\mathcal P})$ be an random locally convex module over $K$ with base $(\Omega, {\mathcal F}, P)$, $x\in E$ and $M$ a ${\mathcal T}_{\varepsilon, \lambda}$--closed $L^0$--convex nonempty subset of $E$. If $x\notin M$, then there is $f\in E^\ast_{\varepsilon, \lambda}$ such that:\\
(1). $(Ref)(x)>\bigvee\{(Ref)(y): y\in M\}~\hbox{on}~(d(x, M)>0)$;\\
(2). $(Ref)(x)=\bigvee\{(Ref)(y): y\in M\}~\hbox{on}~(d(x, M)>0)^c$.\\
\end{proposition}

\begin{remark} When $(E, {\mathcal P})$ is an $RN$ module, $d(x, M)$ is just the random distance from $x$ to $M$. Thus Proposition 4.1 is best possible from the degree that $f$ separates $x$ from $M$.
\end{remark}

In Proposition 4.1, if the condition that ${\tilde I}_A\{x\}\cap {\tilde I}_AM=\emptyset$ for all $A\in {\mathcal F}$ with $P(A)>0$ is also satisfied, then again by (3) of Lemma 3.8 of \cite{TXG-JFA} we have that $(d(x, M)>0)$ has probability 1, which guarantees the separation of $x$ from $M$ by $f$ with probability $1$, so we have Corollary 4.3 below.

\begin{corollary} Let $(E, {\mathcal P})$ be an random locally convex module over $K$ with base $(\Omega, {\mathcal F}, P)$, $x\in E$ and $M$ a ${\mathcal T}_{\varepsilon, \lambda}$--closed $L^0$--convex nonempty subset of $E$. If ${\tilde I}_A\{x\}\cap {\tilde I}_AM=\emptyset$ for all $A\in {\mathcal F}$ with $P(A)>0$, then there is $f\in E^\ast_{\varepsilon, \lambda}$ such that
$(Ref)(x)>\bigvee\{(Ref)(y): y\in M\}~\hbox{on}~\Omega.$

\end{corollary}

In Proposition 4.1, $f$ is asked to belong to $E^\ast_{\varepsilon, \lambda}$, but the future study of random convex analysis requires an $f\in E^\ast_c$ to separate a point from a ${\mathcal T}_{\varepsilon, \lambda}$--closed $L^0$--convex subset, so here we make use of the relation that $E^\ast_{\varepsilon, \lambda}=H_{cc}(E^\ast_c)$ to obtain the following generalization of Proposition 4.1:

\begin{theorem} Let $(E, {\mathcal P})$ be an $RLC$ module over $K$ with base $(\Omega, {\mathcal F}, P)$, $x\in E$ and $M\subset E$ a nonempty ${\mathcal T}_{\varepsilon, \lambda}$--closed $L^0$--convex subset. If $x\notin M$, then there are $f\in E^\ast_c$ and some $B\in {\mathcal F}$ with $P(B)>0$ such that:\\
(1). $(Ref)(x)>\bigvee\{(Ref)(y): y\in M\}~\hbox{on}~B$;\\
(2). $(Ref)(x)=\bigvee\{(Ref)(y): y\in M\}~\hbox{on}~B^c$.
\end{theorem}

\begin{proof}By Proposition 4.1, there exists $g\in E^\ast_{\varepsilon, \lambda}$ such that, by letting $A=(d(x, M)>0)$, \\
(3). $(Reg)(x)>\bigvee\{(Reg)(y): y\in M\}~\hbox{on}~A$;\\
(4). $(Reg)(x)=\bigvee\{(Reg)(y): y\in M\}~\hbox{on}~A^c$.\\

Since $E^\ast_{\varepsilon, \lambda}=H_{cc}(E^\ast_c)$ by Theorem 3.7, $g=\sum^\infty_{n=1}{\tilde I}_{A_n}g_n$ for some countable partition $\{A_n, n\in N\}$ of $\Omega$ to ${\mathcal F}$ and some sequence $\{g_n, n\in N\}$ in $E^\ast_c$. Let $n_0\in N$ be such that $P(A\cap A_{n_0})>0$ and further let $B=A\cap A_{n_0}$ and $f={\tilde I}_{A\cap A_{n_0}}g_{n_0}$, then $f$ and $B$ meet the needs of (1) and (2).
\end{proof}

\subsection{Separation under the locally $L^0$-convex topology}

For the proof of Proposition 4.6 below, let us first recall the following:

\begin{proposition}($See$ \cite{TXG-JFA}). Let $(E, {\mathcal P})$ be an $RLC$ module and $G$ a subset of $E$ such that $G$ has the countable concatenation property. Then ${\bar G}_{\varepsilon, \lambda}={\bar G}_c$, where ${\bar G}_{\varepsilon, \lambda}$ and ${\bar G}_c$ denotes the ${\mathcal T}_{\varepsilon, \lambda}$--  and ${\mathcal T}_c$-- closures of $G$, respectively.
\end{proposition}

By Proposition 4.5 $M$ in Proposition 4.6 below is ${\mathcal T}_{\varepsilon, \lambda}$--closed, further by Corollary 3.6 we have that $E^\ast_c=E^\ast_{\varepsilon, \lambda}$ for $(E, {\mathcal P})$ in Proposition 4.6 below. So we can directly obtain Proposition 4.6 and Corollary 4.7 below from Proposition 4.1 and Corollary 4.3, respectively.

\begin{proposition} Let $(E, {\mathcal P})$ be an $RLC$ module over $K$ with base $(\Omega, {\mathcal F}, P)$ such that ${\mathcal P}$ has the countable concatenation property, $x\in E$ and $M\subset E$ a nonempty ${\mathcal T}_c$--closed $L^0$--convex set with the countable concatenation property. If $x\notin M$, then there is $f\in E^\ast_c$ such that:\\
(1). $(Ref)(x)>\bigvee\{(Ref)(y): y\in M\}~\hbox{on}~(d(x, M)>0)$;\\
(2). $(Ref)(x)=\bigvee\{(Ref)(y): y\in M\}~\hbox{on}~(d(x, M)>0)^c$.\\

\end{proposition}

\begin{corollary} Let $(E, {\mathcal P})$ be an $RLC$ module over $K$ with base $(\Omega, {\mathcal F}, P)$ such that ${\mathcal P}$ has the countable concatenation property, $x\in E$ and $M$ a nonempty ${\mathcal T}_c$--closed $L^0$--convex subset with the countable concatenation property. If ${\tilde I}_A\{x\}\cap {\tilde I}_AM=\emptyset$ for all $A\in {\mathcal F}$ with $P(A)>0$, then there is $f\in E^\ast_c$ such that $$(Re f)(x)>\bigvee\{(Ref)(y): y\in M\} ~\hbox{on}~\Omega.$$
\end{corollary}

Example 4.11 below shows that Corollary 4.7 may be not valid if $M$ lacks the countable concatenation property, so we also correct Theorem 2.8 of \cite{FKV}.

The further study of random convex analysis will need another generalization of Proposition 4.6, namely Corollary 4.8 below, in which the condition that ${\mathcal P}$ has the countable concatenation property will be removed but (1) of Proposition 4.6 will only hold on a subset $B$ of $(d(x, M)>0)$ with $P(B)>0$.

\begin{corollary} Let $(E, {\mathcal P})$ be an $RLC$ module over $K$ with base $(\Omega, {\mathcal F}, P)$, $x\in E$ and $M\subset E$ a nonempty ${\mathcal T}_c$--closed $L^0$--convex set with the countable concatenation property. If $x\notin M$, then there exist $f\in E^\ast_c$ and $B\in {\mathcal F}$ with $P(B)>0$ such that:\\
(1). $(Ref)(x)>\bigvee\{(Ref)(y): y\in M\}~\hbox{on}~B$;\\
(2). $(Ref)(x)=\bigvee\{(Ref)(y): y\in M\}~\hbox{on}~B^c$.
\end{corollary}

\begin{proof} We consider the separation problem in $(E, {\mathcal P}_{cc})$. Since ${\mathcal P}_{cc}$ has the countable concatenation property and the locally $L^0$--convex topology induced by ${\mathcal P}_{cc}$ is stronger than that induced by ${\mathcal P}$. We can apply Proposition 4.6 to $(E, {\mathcal P}_{cc})$, $x$ and $M$, then there is $g\in (E, {\mathcal P}_{cc})^\ast_c$ such that, by letting $A=(d(x, M)>0)$,\\
(3). $(Reg)(x)>\bigvee\{(Reg)(y): y\in M\}~\hbox{on}~A$;\\
(4). $(Reg)(x)=\bigvee\{(Reg)(y): y\in M\}~\hbox{on}~A^c$.\\
Here, please note that ${\mathcal P}$ and ${\mathcal P}_{cc}$ induce the same $d(x, M)$, so $A$ is still a representative of $[d(x, M)>0]$.

Since $(E, {\mathcal P}_{cc})^\ast_c=H_{cc}(E^\ast_c)$ by Theorem 3.7, $g=\sum^\infty_{n=1}{\tilde I}_{A_n}g_n$ for some countable partition $\{A_n, n\in N\}$ of $\Omega$ to ${\mathcal F}$ and some sequence $\{g_n, n\in N\}$ in $E^\ast_c$. Let $n_0\in N$ be such that $P(A\cap A_{n_0})>0$ and further let $B=A\cap A_{n_0}$ and $f={\tilde I}_{A\cap A_{n_0}}g_{n_0}$, then $f$ and $B$ meet the needs of (1) and (2).
\end{proof}

\begin{remark} Let $\xi$ be any element in $L^0({\mathcal F}, K)$ and $\xi_0$ a representative of $\xi$. Define $\xi_0^{-1}\colon\Omega\to K$ by $\xi_0^{-1}(\omega)=(\xi_0(\omega))^{-1}$ if $\xi_0(\omega)\neq 0$ and by $0$ if $\xi_0(\omega)=0$, then $\xi^{-1}:=$ the equivalence class of $\xi_0^{-1}$ is called the generalized inverse of $\xi$. $|\xi|^{-1}\xi$ is called the sign of $\xi$, denoted by $sgn(\xi)$, then ${\overline {sgn(\xi)}}\xi=|\xi|$, where ${\overline {sgn(\xi)}}$ stands for the complex conjugate of $sgn(\xi)$. Further, we also have that $\xi\cdot\xi^{-1}=\xi^{-1}\cdot\xi=I_{[\xi\neq 0]}$. If $M$ in Theorem 4.4 or Corollary 4.8 is also $L^0$--balanced, then one can make use of the notion of the sign for elements in $L^0({\mathcal F}, K)$ to see that (1) and (2) of the two results can be rewritten as ( cf. \cite{TXG-strict} ):\\
(1). $|f(x)|>\bigvee\{|f(y)|: y\in M\}~\hbox{on}~B$;\\
(2). $|f(x)|=\bigvee\{|f(y)|: y\in M\}~\hbox{on}~B^c$.

Let $\xi=|f(x)|$ and $\eta=\bigvee\{|f(y)|: y\in M\}$, then multiplying the above two sides by $(\frac{\xi+\eta}{2})^{-1}$ and replacing $f$ with $(\frac{\xi+\eta}{2})^{-1}f$ ( still denoted by $f$ ) will obtain the following two relations:\\
(3). $|f(x)|>\bigvee\{|f(y)|: y\in M\}$;\\
(4). $|f(x)|\nleqslant 1$ and $\bigvee\{|f(y)|: y\in M\}\leq 1$.\\
(3) and (4) will be used in the proof of random bipolar theorem in our forthcoming paper.

\end{remark}

To study the properties of proper $L^0$-convex lower semicontinuous functions on an $RLC$ module under the locally $L^0$-convex topology, we need the following separation proposition by neighborhoods between a point and a $\mathcal{T}_c$-closed set, namely Theorem 4.10 below.

\begin{theorem} Let $(E,\mathcal{P})$ be an $RLC$ module over $K$ with base $(\Omega,\mathcal{F},P)$ such that $\mathcal{P}$ has the countable concatenation property, $M$ a $\mathcal{T}_c$--closed subset with the countable concatenation property  and $x\in E$ such that $\tilde{I}_A\{x\}\bigcap\tilde{I}_A M=\emptyset$ for all $A\in \mathcal{F}$ with $P(A)>0$. Then there is an $L^0$--convex, $L^0$--absorbent and $L^0$--balanced $\mathcal{T}_c$--neighborhood $U$ of $\theta$ such that $$\tilde{I}_A(x+U)\bigcap\tilde{I}_A(M+U)=\emptyset$$ for all $A\in\mathcal{F}$ with $P(A)>0$.
\end{theorem}

Example 4.11 below shows that Theorem 4.10 may be not valid if $M$ is only required to satisfy the condition that $\tilde{I}_{A}M+\tilde{I}_{A^c}M\subset M$ for all $A\in\mathcal{F}$, so we also correct Lemma 2.28 of \cite{FKV}.

To prove Theorem 4.10, let us first recall Lemma 3.10 of \cite{TXG-JFA}: let $(E,\mathcal{P})$ be an $RLC$ module over $K$ with base $(\Omega,\mathcal{F},P)$, for each $\mathcal{Q}\in\mathcal{P}_f$ and $\varepsilon\in L^{0}_{++}(\mathcal{F})$, let $U_{\mathcal{Q},\varepsilon}[x]=\{y\in E~|~\|x-y\|_{\mathcal{Q}}\leq\varepsilon\}$, $e_{\mathcal{Q}}(x,M)=\bigwedge\{\varepsilon\in L^0_{++}(\mathcal{F})~|~U_{\mathcal{Q},\varepsilon}[x]\bigcap M\neq\emptyset\}$ and $e(x,M)=\bigvee\{e_{\mathcal{Q}}(x,M)~|~\mathcal{Q}\in\mathcal{P}_{f}\}$. Then $d(x,M)=e(x,M)$ and further $e(x,M)\bigwedge 1\in L^0_{++}(\mathcal{F})$ for $x$ and $M$ in Theorem 4.10.

We can now prove Theorem 4.10.

\noindent{\bf Proof of Theorem 4.10.}
We can assume that $M\neq\emptyset$ and $x=\theta$ (otherwise by translation), it suffices to construct an $L^{0}$-convex, $L^{0}$-absorbent and $L^{0}$-balanced neighborhood $U$ of $\theta$ such that $$\tilde{I}_{A}U\bigcap\tilde{I}_{A}(M+U)=\emptyset$$ for all $A\in \mathcal{F}$ with $P(A)>0$. Let $\varepsilon^{\ast}=1\wedge e(\theta,M)$, then $\varepsilon^{\ast}\in L^{0}_{++}$ by Lemma 3.10 of \cite{TXG-JFA}. We will show that $\varepsilon^{\ast}$ satisfies:\\
(i) There is an $L^0$-seminorm $\|\cdot\|^{\ast}\in\mathcal{P}$ such that $\frac{\varepsilon^{\ast}}{2}<\bigwedge\{\varepsilon\in L^{0}_{++}~|~N_{\theta}(\|\cdot\|^{\ast},\varepsilon)\bigcap M\neq\emptyset\}~ on~\Omega$, where $N_{\theta}(\|\cdot\|^{\ast},\varepsilon)=\{y\in E~|~\|y\|^{\ast}\leq\varepsilon\}$.\\
(ii) $\tilde{I}_{A}N_{\theta}(\|\cdot\|^{\ast},\frac{\varepsilon^{\ast}}{2})\bigcap\tilde{I}_{A}M=\emptyset$ for all $A\in \mathcal{F}$. (Note that $N_{\theta}(\|\cdot\|^{\ast},\frac{\varepsilon^{\ast}}{2})$ is $L^{0}$-convex, $L^{0}$-absorbent, $L^{0}$-balanced and $\mathcal{T}_c$-closed.)

To prove (i), for all finite $\mathcal{Q}\in \mathcal{P}$, let $\varepsilon_{\mathcal{Q}}=\bigwedge\{\varepsilon\in L^{0}_{++}~|~N_{\theta}(\mathcal{Q},\varepsilon)\bigcap M\neq\emptyset\}$, where $N_{\theta}(\mathcal{Q},\varepsilon)=\{y\in E~|~\|y\|_{\mathcal{Q}}\leq\varepsilon\}$. For finite $\mathcal{Q}$, $\mathcal{Q}^{\prime}\in \mathcal{P}$, $N_{\theta}(\mathcal{Q}\bigcup\mathcal{Q}^{\prime},\varepsilon)\subset N_{\theta}(\mathcal{Q},\varepsilon),N_{\theta}(\mathcal{Q}^{\prime},\varepsilon)$. Thus, the collection $\{\varepsilon_{\mathcal{Q}}~|~\mathcal{Q}\subset\mathcal{P}$ finite $\}$ is directed upwards and hence there is an increasing sequence $(\varepsilon_{\mathcal{Q}_{n}})$ with $1\wedge\varepsilon_{\mathcal{Q}_{n}}\nearrow \varepsilon^{\ast}$. Let $B_n$ be a representative of $[\varepsilon_{\mathcal{Q}_{n}}>\frac{\varepsilon^{\ast}}{2}]$ for any $n\in N$, $A_1=B_1$ and $A_n=B_n\setminus B_{n-1}$ for any $n\geq 2$. Then we can, without loss of generality, assume $\bigcup_{n\in N}A_n=\Omega$ since $\varepsilon^{\ast}>\frac{\varepsilon^{\ast}}{2}$ on $\Omega$. Further, the $L^0$-seminorm $\|\cdot\|^{\ast}=\sum\limits_{n\in N}\tilde{I}_{A_n}\|\cdot\|_{\mathcal{Q}_n}$ is an element of $\mathcal{P}$ since $\mathcal{P}$ has the countable concatenation property.

Finally, to prove (ii), assume there is $A\in \mathcal{F}$, $P(A)>0$ and $y\in M$ such that $\tilde{I}_{A}y\in \tilde{I}_{A}N_{\theta}(\|\cdot\|^{\ast},\frac{\varepsilon^{\ast}}{2})$, then $\tilde{I}_{A}\bigwedge\{\varepsilon\in L^{0}_{++}~|~N_{\theta}(\|\cdot\|^{\ast},\varepsilon)\bigcap M\neq\emptyset\}\leq\tilde{I}_{A}\frac{\varepsilon^{\ast}}{2}$ in contradiction to the statement in (i).

To sum up, from the above proofs we have $\|\cdot\|\in\mathcal{P}$ and $\varepsilon\in L^{0}_{++}$ such that $\tilde{I}_{A}N_{\theta}(\|\cdot\|,\varepsilon)\bigcap\tilde{I}_{A}M=\emptyset$ for all $A\in\mathcal{F}$ with $P(A)>0$. This implies $\tilde{I}_{A}N_{\theta}(\|\cdot\|,\frac{\varepsilon}{2})\bigcap\tilde{I}_{A}(M+N_{\theta}(\|\cdot\|,\frac{\varepsilon}{2})=\emptyset$ for all $A\in\mathcal{F}$ with $P(A)>0$ and the assertion follows.  \hfill $\square$

\begin{example} Let $(\Omega,\mathcal{F},P)$ be a nonatomic probability space (namely $\mathcal{F}$ does not include any $P$--atoms), $(E,\mathcal{P})=(L^{0}(\mathcal{F},R),|\cdot|)$ and $M=\{x\in E~|~$there exists a positive number $m_x$ such that $x>m_x$ on $\Omega\}$. Then Claim 4.12 below shows that $M$ is $L^0$--convex, $\mathcal{T}_c$--closed and $\mathcal{T}_c$--open. Further, Claim 4.13 below shows that $\tilde{I}_A\{0\}\bigcap\tilde{I}_A M=\emptyset$ for all $A\in\mathcal{F}$ with $P(A)>0$, but for each $L^0$--convex, $L^0$--absorbent and $L^0$--balanced $\mathcal{T}_c$--neighborhood $U$ of 0 there is an $A_U\in\mathcal{F}$ with $P(A_U)>0$ such that $$\tilde{I}_{A_U}U\bigcap\tilde{I}_{A_U}(M+U)\neq\emptyset.$$
\end{example}

\begin{claim} $M$ in Example 4.11 is $L^0$--convex, $\mathcal{T}_c$--closed and $\mathcal{T}_c$--open.
\end{claim}

\begin{proof} First, it is obvious that $M$ is $L^0$--convex.

Second, $M$ is $\mathcal{T}_c$--open. For any $y\in M$, by definition there is some positive number $m_y$ such that $y>m_y$ on $\Omega$. Let $\varepsilon^0\equiv\frac{1}{2}m_y$ and $\varepsilon$ be the equivalence class of $\varepsilon^0$, then $\varepsilon\in L^0_{++}(\mathcal{F})$ and hence $B(\varepsilon):=\{x\in E~|~|x|\leq\varepsilon\}$ is a $\mathcal{T}_c$--neighborhood of 0, it is also easy to check that $y+B(\varepsilon)\subset M$.

Finally, $M$ is also $\mathcal{T}_c$--closed, namely $E\setminus M$ is $\mathcal{T}_c$--open, which will be proved in the following three cases.

Case (1): when $y\in E\setminus M$ and $y\not\in L^0_+(\mathcal{F})$, there is $D\in\mathcal{F}$ with $P(D)>0$ such that $y<0$ on $D$. Let $\varepsilon=\tilde{I}_{D^c}+\frac{1}{2}\tilde{I}_{D}|y|(\in L^{0}_{++}(\mathcal{F}))$ and $B(\varepsilon)=\{x\in E~|~|x|\leq\varepsilon\}$, then $y+B(\varepsilon)\subset E\setminus M$. In fact, for any $z\in y+B(\varepsilon)$, $z-y\leq \tilde{I}_{D^c}+\frac{1}{2}\tilde{I}_{D}|y|$ implies that $z\leq y+\frac{1}{2}|y|=-\frac{1}{2}|y|<0$ on $D$, namely $z\in E\setminus M$.

Case (2): when $y\in E\setminus M$, $y\in L^0_+(\mathcal{F})$ and $y\not\in L^{0}_{++}(\mathcal{F})$, then there is $D\in\mathcal{F}$ with $P(D)>0$ such that $y=0$ on $D$. Since $(\Omega,\mathcal{F},P)$ is nonatomic, there is a countable partition $\{D_n,n\in N\}$ of $D$ to $\mathcal{F}$ such that $P(D_n)=\frac{1}{2^n}P(D)$ for each $n\in N$. Let $\varepsilon=\tilde{I}_{D^c}+\Sigma_{n=1}^{\infty}\frac{1}{n}\tilde{I}_{D_n} (\in L^{0}_{++}(\mathcal{F}))$ and $B(\varepsilon)=\{x\in E~|~|x|\leq\varepsilon\}$, then $z\leq \frac{1}{n}$ on $D_n$ for any $z\in y+B(\varepsilon)$, which implies that $P\{\omega\in\Omega~|~z(\omega)\leq\frac{1}{n}\}\geq P(D_n)>0$ for all $n\in N$, namely $y+B(\varepsilon)\subset E\setminus M$.

Case (3): when $y\in E\setminus M$ and $y\in L^0_{++}(\mathcal{F})$, then $P\{\omega\in \Omega~|~y(\omega)<\frac{1}{n}\}>0$ for each $n\in N$ by the definition of $M$. Let $H_n=[y<\frac{1}{n}]$ and $D_n=[\frac{1}{n+1}\leq y<\frac{1}{n}]$ for any $n\in N$, then $D_i\bigcap D_j=\emptyset$ for $i\neq j$ and $H_n=\Sigma_{i=n}^{\infty}D_i$. Obviously, it is impossible that there is some $k\in N$ such that $P(D_n)=0$ for all $n\geq k$. So, we can suppose, without loss of generality, that $P(D_n)>0$ for each $n\in N$. Let $D=\Sigma_{n=1}^{\infty}D_n$, $\varepsilon=I_{D^{c}}+\Sigma_{n=1}^{\infty} \frac{1}{n}I_{D_{n}}(\in L^{0}_{++}(\mathcal{F}))$ and $B(\varepsilon)
=\{x\in E|~|x|\leq \varepsilon\}$, then for any $z\in y+B(\varepsilon)$, $z\leq\frac{2}{n}$ on $D_n$, which means that $P\{\omega\in \Omega|z(\omega)\leq\frac{2}{n}\}\geq P(D_n)>0$ for each $n\in N$, and hence $z\in E\setminus M$.
\end{proof}

\begin{claim} Let $(E, {\mathcal P})$ and $M$ be the same as in Example 4.11. Then ${\tilde I}_A\{0\}\cap {\tilde I}_AM=\emptyset$ for all $A\in {\mathcal F}$ with $P(A)>0$. But for any $L^0$--convex, $L^0$--absorbent and $L^0$--balanced ${\mathcal T}_c$--neighborhood $U$ of $0$ there is always $A_U\in {\mathcal F}$ with $P(A_U)>0$ such that ${\tilde I}_{A_U}U\cap {\tilde I}_{A_U}(M+U)\neq\emptyset$.
\end{claim}

\begin{proof} There is $\varepsilon \in L^0_{++}({\mathcal F})$ for $U$ stated above such that $B(\varepsilon):=\{x\in E~|~|x|\leq\varepsilon\}\subset U$. For a representative $\varepsilon^0$ of $\varepsilon$, let $A_1=\{\omega\in \Omega~|~\varepsilon^0(\omega)\geq 1\}$ and $A_n=\{\omega\in\Omega~|~\frac{1}{n}\leq \varepsilon^0(\omega)<\frac{1}{n-1}\}$ for $n\geq 2$, then it is clear that $\sum^{\infty}_{n=1}P(A_n)=1$, and hence there is some $n_0\in N$ such that $P(A_{n_0})>0$. Let $A_U=A_{n_0}$ and $y_0={\tilde I}_{A^c_U}+{\tilde I}_{A_U}\varepsilon$, then $\frac{1}{n_0}\leq y_0<\frac{1}{n_0-1}$ on $A_U$ ( note: this is also true for $n_0=1$ ) and $y_0\geq \frac{1}{n_0}$ on $\Omega$ (namely, $y_0\in M$). Since ${\tilde I}_{A_U}y_0={\tilde I}_{A_U}\varepsilon\in {\tilde I}_{A_U}B(\varepsilon)\subset {\tilde I}_{A_U} U$ and ${\tilde I}_{A_U}y_0\in {\tilde I}_{A_U}M\subset{\tilde I}_{A_U}(M+U)$, so ${\tilde I}_{A_U}U\cap {\tilde I}_{A_U}(M+U)\neq\emptyset$.
\end{proof}

Example 4.11 also shows that Corollary 4.7 may be not valid if $M$ lacks the countable concatenation property. Since $(E, {\mathcal P})=(L^0({\mathcal F}, R), |\cdot|)$ is an $RN$ module, $|\cdot|$ has the countable concatenation property and $E^\ast_c=E^\ast_{\varepsilon, \lambda}$. It is obvious that $0\in {\overline{M}}_{\varepsilon, \lambda}$ ( namely, the ${\mathcal T}_{\varepsilon, \lambda}$--closure of $M$ ), and hence for each $f\in E^\ast_c=E^\ast_{\varepsilon, \lambda}$ there exists a sequence $\{y_n, n\in N\}$ in $M$ such that $\{f(y_n): n\in N\}$ converges in probability P to $0$, which means that it is impossible that there exists $f\in E^\ast_c$ such that $0=f(0)>\bigvee\{f(y): y\in M\}$ on $\Omega$.

\section{The Fenchel-Moreau dual representation theorems in random locally convex modules under the two kinds of topologies}

Let $E$ be an $L^0({\mathcal F})$--module and $f$ a function from $E$ to ${\bar L}^0({\mathcal F})$. The effective domain of $f$ is denoted by $dom(f):=\{x\in E~|~f(x)<+\infty~\hbox{on}~\Omega\}$ and the epigraph of $f$ by $epi(f):=\{(x,r)\in E\times L^0({\mathcal F})~|~f(x)\leq r\}$. $f$ is proper if $dom(f)\neq\emptyset$ and $f(x)>-\infty~\hbox{on}~\Omega$. $f$ is $L^0$--convex if $f(\xi x+(1-\xi)y)\leq \xi f(x)+(1-\xi)f(y)$ for all $x,~y\in E$ and $\xi\in L^0_+({\mathcal F})$ with $0\leq \xi\leq 1$, where the following convention is adopted: $0\cdot(\pm\infty)=0$ and $+\infty\pm(\pm\infty)=+\infty$. $f\colon E\to {\bar L}^0({\mathcal F})$ is said to be local ( or, to have the local property ) if ${\tilde I}_Af(x)={\tilde I}_Af({\tilde I}_Ax)$ for all $x\in E$ and $A\in {\mathcal F}$. In \cite{FKV-appro}, it is proved that an $L^0$-convex function is local.

The main results of this section are Theorems 5.3 and 5.5 below, which are the Fenchel-Moreau dual representation theorems in random locally convex modules under the two kinds topologies, respectively.

\subsection{The Fenchel-Moreau dual representation theorem in random locally convex modules under the $(\varepsilon,\lambda)$--topology}

\begin{lemma}($See$ \cite{FKV,FKV-appro}). Let $E$ be an $L^0({\mathcal F})$--module. Then a proper function $f\colon E\to {\bar L}^0({\mathcal F})$ is $L^0$--convex iff $f$ is local and $epi(f)$ is $L^0$--convex.
\end{lemma}

\begin{definition} Let $(E,\mathcal{P})$ be an $RLC$ module over $R$ with base $(\Omega,\mathcal{F},P)$ and $f:E\rightarrow \bar{L}^{0}(\mathcal{F})$ a proper $L^{0}$--convex function. $f$ is $\mathcal{T}_{\varepsilon,\lambda}$--lower semicontinuous if $epi(f)$ is closed in $(E,\mathcal{T}_{\varepsilon,\lambda})\times(L^0(\mathcal{F}),\mathcal{T}_{\varepsilon,\lambda})$.
\end{definition}

As usual, let $(E,\mathcal{P})$ be an $RLC$ module over $R$ with base $(\Omega,\mathcal{F},P)$ and $f:E\rightarrow L^0(\mathcal{F})$, $f$ is $\mathcal{T}_{\varepsilon,\lambda}$--continuous if $f$ is continuous from $(E,\mathcal{T}_{\varepsilon,\lambda})$ to $(L^0(\mathcal{F}),\mathcal{T}_{\varepsilon,\lambda})$.

As to why we adopt Definition 5.2 for the $\mathcal{T}_{\varepsilon,\lambda}$--lower semicontinuity of an $L^0$--convex function, we interpret this as follows. If we define the $\mathcal{T}_{\varepsilon,\lambda}$--lower semicontinuity of a proper function $f:(E,\mathcal{P})\rightarrow\bar{L}^0({\mathcal{F}})$ via ``$\{x\in E~|~f(x)\leq r\}$ is $\mathcal{T}_{\varepsilon,\lambda}$--closed for all $r\in L^0(\mathcal{F})$", then this notion is too weak to meet some natural needs of other topics as in \cite{TXG-YJY}. If we define $f$ to be lower semicontinuous via ``$\underline{lim}_\alpha f(x_\alpha):=\bigvee_{\beta\in\Gamma}(\bigwedge_{\alpha\geq\beta}f(x_{\alpha}))\geq f(x)$ for all nets $\{x_{\alpha},\alpha\in\Lambda\}$ in $E$ such that it converges in the $(\varepsilon,\lambda)$--topology to some $x\in E$", the notion is, however, meaningless in the random setting, since we can construct a real $RLC$ module $(E,\mathcal{P})$ and a $\mathcal{T}_{\varepsilon,\lambda}$--continuous $L^0$--convex function $f$ from $E$ to $L^0(\mathcal{F})$, whereas $f$ is not a lower semicontinuous function under this notion. In fact, Definition 5.2 has been proved natural and fruitful, see \cite{TXG-YJY} or this subsection.

The proof of Theorem 5.3 below, namely the random version under the $(\varepsilon,\lambda)$--topology of the classical Fenchel-Moreau duality theorem, is more complicated since the complicated stratification structure in the random setting needs to be considered. As compared with Theorem 5.5 below, namely the random version under the locally $L^0$--topology of the classical Fenchel-Moreau duality theorem, Theorem 5.3 is more natural since it
has the same shape as the classical Fenchel-Moreau duality theorem.

Let $(E,\mathcal{P})$ be an $RLC$ module over $R$ with base $(\Omega,\mathcal{F},P)$ and $f:E\rightarrow \bar{L}^0(\mathcal{F})$ a proper $\mathcal{T}_{\varepsilon,\lambda}$--lower semicontinuous $L^0$--convex function. We define $f^{\ast}_{\varepsilon,\lambda}:E^{\ast}_{\varepsilon,\lambda}\rightarrow\bar{L}^{0}(\mathcal{F})$ by $f^{\ast}_{\varepsilon,\lambda}(g)=\bigvee\{g(x)-f(x)~|~x\in E\}$ for all $g\in E^{\ast}_{\varepsilon,\lambda}$, called the $\mathcal{T}_{\varepsilon,\lambda}$-conjugate function of $f$, and $f^{\ast\ast}_{\varepsilon,\lambda}:E\rightarrow\bar{L}^0(\mathcal{F})$ by $f^{\ast\ast}_{\varepsilon,\lambda}(x)=\bigvee\{g(x)-f^{\ast}_{\varepsilon,\lambda}(g)~|~g\in E^{\ast}_{\varepsilon,\lambda}\}$ for all $x\in E$, called the $\mathcal{T}_{\varepsilon,\lambda}$-bi-conjugate function of $f$.

\begin{theorem} Let $(E,\mathcal{P})$ be an $RLC$ module over $R$ with base $(\Omega,\mathcal{F},P)$ and $f:E\rightarrow \bar{L}^0(\mathcal{F})$ a proper $\mathcal{T}_{\varepsilon,\lambda}$--lower semicontinuous $L^0$--convex function. Then $f^{\ast\ast}_{\varepsilon,\lambda}=f$.
\end{theorem}

\begin{proof} We fix $x_0\in E$ and proceed in two steps below.

Step 1: Let $\beta\in L^{0}(\mathcal{F})$ with $\beta<f(x_0)$ on $\Omega$. In this step, we show there is a continuous function $h:(E,\mathcal{T}_{\varepsilon,\lambda})\rightarrow (L^0(\mathcal{F}),\mathcal{T}_{\varepsilon,\lambda})$ of the form $$h(x)=g(x)+z,$$ where $g\in E^{\ast}_{\varepsilon,\lambda}$ and $z\in L^{0}(\mathcal{F})$, such that $h(x_0)=\beta$ and $h(x)\leq f(x)$ for all $x\in E$. To this end, we separate $(x_0,\beta)$ from $epi(f)$ by means of Corollary 4.3. It applies since $\beta<f(x_0)$ on $\Omega$ and the local property of $f$ imply $\tilde{I}_A(x_0,\beta)\bigcap \tilde{I}_A epi(f)=\emptyset$ for all $A\in \mathcal{F}$ with $P(A)>0$. (Note, $epi(f)$ is closed in $(E,\mathcal{T}_{\varepsilon,\lambda})\times(L^0(\mathcal{F}),\mathcal{T}_{\varepsilon,\lambda})$ by Definition 5.2.) Hence, there are $g_1\in E^{\ast}_{\varepsilon,\lambda}$ and $g_2\in (L^{0}(\mathcal{F}))^{\ast}_{\varepsilon,\lambda}$ ( in fact, $(L^{0}(\mathcal{F}))^{\ast}_{\varepsilon,\lambda}=L^{0}(\mathcal{F})$) such that $\delta=\bigvee\limits_{(x,y)\in epi(f)}(g_1(x)+g_2(y))<g_1(x_0)+g_2(\beta)$ on $\Omega$. This has two consequences:\\
(i) $g_2(1)\leq 0$.

Indeed, $g_2(y)=yg_2(1)$ for all $y\in L^{0}(\mathcal{F})$. Further, $(x,y)\in epi(f)$ for arbitrary large $y$ as long as $f(x)\leq y$. Hence, for large $y\in L^{0}(\mathcal{F})$, $g_1(x)+g_2(y)$ is large on $[g_2(1)>0]$ and yet bounded above by $g_1(x_0)+g_2(\beta)$. This implies $P([g_2(1)>0])=0$.\\
(ii) $[f(x_0)<+\infty]\subset [g_2(1)<0]$.

Indeed, define $\tilde{x}_0=I_{[f(x_0)<+\infty]}x_0+I_{[f(x_0)=+\infty]}x$ for some $x\in dom(f)$. By $L^0$-convexity of $f$, $\tilde{x}_0\in dom(f)$. Local property of $f$ and the definition of $\delta$ imply on $[f(x_0)<+\infty]$ $$g_1(x_0)+g_2(f(x_0))=g_1(\tilde{x}_0)+g_2(f(\tilde{x}_0))<g_1(x_0)+g_2(\beta).$$ Hence, $f(x_0)g_2(1)=g_2(f(x_0))<g_2(\beta)=\beta g_2(1)$ on $[f(x_0)<+\infty]$ and so $g_2(1)<0$ on $[f(x_0)<+\infty]$.

We distinguish the cases $x_0\in dom(f)$ and $x_0\not\in dom(f)$.

Case 1. Assume $x_0\in dom(f)$. By (ii), $g_2(1)<0$ on $\Omega$. Thus, define $h$ by $$h(x)=-\frac{g_1(x-x_0)}{g_2(1)}+\beta$$ for all $x\in E$, which is as required. Indeed, $h(x)\leq f(x)$ for all $x\in dom(f)$ by the definition of $\delta$. If $x\not\in dom(f)$ we have $$I_Bh(x)=I_Bh(x^{\prime})\leq I_Bf(x^{\prime})=I_Bf(x),$$ where $x^{\prime}=I_Bx+I_{B^{c}}x^{\prime\prime}$ for some $x^{\prime\prime}\in dom(f)$ and $B=[f(x)<+\infty]$. Hence, $h(x)\leq f(x)$ for all $x\in E$.

Case 2. Assume $x_0\not\in dom(f)$. Then choose any $x_{0}^{\prime}\in dom(f)$ and let $h^{\prime}$ be the corresponding $L^{0}(\mathcal{F})$-affine minorant as constructed in case 1 above. Define $A_1=[g_2(1)<0]$, $A_2=A_{1}^{c}$ and $h_1,h_2:E\rightarrow L^{0}(\mathcal{F})$,$$h_1(x)=I_{A_1}(-\frac{g_1(x-x_0)}{g_2(1)}+\beta),$$

$$h_2(x)=I_{A_2}[h^{\prime}(x)+I_{[h^{\prime}(x_0)\geq \beta]}(\beta-h^{\prime}(x_0))+I_{[h^{\prime}(x_0)<\beta]}\frac{\beta-h^{\prime}(x_0)}{\tilde{h}(x_0)}\tilde{h}(x)]$$
with the convention $\frac{0}{0}=0$, where $\tilde{h}:E\rightarrow L^0(\mathcal{F})$ is defined by $$\tilde{h}(x)=\delta-g_1(x).$$ Note that on $[g_2(1)=0]$, $\tilde{h}(x_0)<0$ and $\tilde{h}(x)\geq 0$ for all $x\in dom(f)$. It follows that $$h=h_1+h_2$$ is as required.

Step 2: It is clear that $f\geq f_{\varepsilon,\lambda}^{\ast\ast}$ by the definition of $f_{\varepsilon,\lambda}^{\ast\ast}$. By way of contradiction, assume $f(x_0)>f_{\varepsilon,\lambda}^{\ast\ast}(x_0)$ on a set of positive measure. Then there is $\beta\in L^{0}(\mathcal{F})$ with $\beta>f_{\varepsilon,\lambda}^{\ast\ast}(x_0)$ on a set of positive measure and $\beta<f(x_0)$ on $\Omega$. The first step of this proof yields $h:E\rightarrow L^{0}(\mathcal{F})$, $$h(x)=g(x)+z$$ for all $x\in E$, for $g\in E^{\ast}_{\varepsilon,\lambda}$ and $z\in L^{0}(\mathcal{F})$, such that $h(x_0)=\beta$ and $h(x)\leq f(x)$ for all $x\in E$. We derive a contradiction as $$f_{\varepsilon,\lambda}^{\ast\ast}(x_0)\geq g(x_0)-f_{\varepsilon,\lambda}^{\ast}(g)$$ $$~~~~~~~~~~~~~~~~~~~~~~=g(x_0)-\bigvee_{x\in E}(g(x)-f(x))$$ $$~~~~~~~~~~~~~~~~~~~~~~~~~~~~\geq g(x_0)-\bigvee_{x\in E}(g(x)-h(x))=\beta.$$
\end{proof}

\begin{remark}  Historically, the random generalization of the classical Fenchel-Moreau duality theorem is first studied under the framework of a locally $L^0$--convex module in \cite{FKV}, where Theorem 3.8 of \cite{FKV} was given and some good contributions were made, for example, the construction of $h$ in the proof of Theorem 5.3 is just taken from the process of the proof of Theorem 3.8 in \cite{FKV}. But, Theorem 2.4 and Lemma 2.28 of \cite{FKV} is not true which makes Theorem 3.8 in \cite{FKV} also wrong. Theorem 5.5 below of this paper is given in order to correct Theorem 3.8 of \cite{FKV}.

\end{remark}

\subsection{The Fenchel-Moreau dual representation theorem in random locally convex modules under the locally $L^{0}$--convex topology}

If $(E, {\mathcal T})$ is a topological $L^0({\mathcal F})$--module, in \cite{FKV} a proper function $f:E\rightarrow \bar{L}^{0}(\mathcal{F})$ is lower semicontinuous (or ${\mathcal T}$--lower semicontinuous if there is a possible confusion) if $\{x\in E~|~f(x)\leq r\}$ is closed for all $r\in L^0({\mathcal F})$. But up to now, we have not seen a strict proof that this kind of ${\mathcal T}$--lower semicontinuity implies the epigraph of $f$ is closed in $(E, {\mathcal T})\times (L^0({\mathcal F}), {\mathcal T}_c)$. For this, in this paper we say that $f$ is ${\mathcal T}$--lower semicontinuous if the epigraph of $f$ is closed in $(E, {\mathcal T})\times (L^0({\mathcal F}), {\mathcal T}_c)$.

We can now state the main result of this subsection as Theorem 5.5 below, which is a modification and improvement of Theorem 3.8 of \cite{FKV}. Let $(E, {\mathcal P})$ be an $RLC$ module over $R$ with base $(\Omega, {\mathcal F}, P)$ such that $E$ has the countable concatenation property. If $f$ is a proper, ${\mathcal T}_c$--lower semicontinuous $L^0({\mathcal F})$--convex function from $E$ to ${\bar L}^0({\mathcal F})$. We define $f^\ast_c\colon E^\ast_c\to {\bar L}^0({\mathcal F})$ by $f^\ast_c(g)=\bigvee\{g(x)-f(x)~|~x\in E\}$ for all $g\in E^\ast_c$, called the ${\mathcal T}_c$--conjugate ( or penalty ) function of $f$, and $f^{\ast\ast}_{c}\colon E \to {\bar L}^0({\mathcal F})$ by $f^{\ast\ast}_{c}(x)=\bigvee\{g(x)-f^{\ast}_{c}(g)~|~g\in E^\ast_c\}$ for all $x\in E$, called the ${\mathcal T}_c$--bi-conjugate function of $f$.

\begin{theorem} Let $(E, {\mathcal P})$ be an $RLC$ module over $R$ with base $(\Omega, {\mathcal F}, P)$ such that $E$ has the countable concatenation property. If $f$ is a proper, ${\mathcal T}_c$--lower semicontinuous $L^0({\mathcal F})$--convex function from $E$ to ${\bar L}^0({\mathcal F})$, then $f^{\ast\ast}_{c}=f$.
\end{theorem}

As compared with Theorem 3.8 of \cite{FKV}, besides, a locally $L^0$--convex module is replaced by a random locally convex module, Theorem 5.5 also requires the additional condition that $E$ has the countable concatenation property and removes the condition that ${\mathcal P}$ has the countable concatenation property. Since the original proof of Theorem 3.8 of \cite{FKV} essentially depends on Corollary 4.7, one can immediately see that the additional condition is essential, whereas Theorem 3.7 can be used to remove the condition on ${\mathcal P}$.
Besides, we remind the reader of the essential distinction between `` the hypothesis that $E$ has the countable concatenation property in our Theorem 5.5 '' and `` the hypothesis that $E$ has the countable concatenation property in Theorem 3.8 of \cite{FKV} ''. According to our Definition 2.14, the hypothesis in our Theorem 5.5 is purely algebraic, whereas, according to Definition 2.7 of \cite{FKV}, the hypothesis in Theorem 3.8 of \cite{FKV} is relative to topology. In fact, Definition 2.7 of \cite{FKV} is not well defined since \cite{FKV} did not give a reasonable interpretation of ``$\sum_{n\in N}I_AU_n$'' in Definition 2.7 of \cite{FKV}.

To prove Theorem 5.5, let us first study the properties of an $L^0$--convex function.

By Proposition 4.5, one can easily see the following:

\begin{lemma} Let $(E, {\mathcal P})$ be an $RLC$ module over $R$ with base $(\Omega, {\mathcal F}, P)$ such that both $E$ and ${\mathcal P}$ have the countable concatenation property. If $f\colon E\to {\bar L}^0({\mathcal F})$ is a proper and local function, then the following are equivalent:\\
(1). $f$ is ${\mathcal T}_c$--lower semicontinuous.\\
(2). $f$ is ${\mathcal T}_{\varepsilon,\lambda}$--lower semicontinuous.\\
\end{lemma}

To prove Theorem 5.5, we still need Lemma 5.7 below, which is almost obvious but frequently used in the proofs of the forthcoming study of the relations among conditional risk measures, and thus we summarize and prove it as follows:

\begin{lemma} Let $E$ be an $L^{0}(\mathcal{F})$--module with the countable concatenation property. Then we have the following statements:

\noindent $(1)$. Let $f:E\rightarrow \bar{L}^{0}(\mathcal{F})$ have the local property and $x=\Sigma_{n=1}^{\infty}\tilde{I}_{A_n}x_n$ for some countable partition $\{A_n,n\in N\}$ of $\Omega$ to $\mathcal{F}$ and some sequence $\{x_n,n\in N\}$ in $E$, then $f(x)=\Sigma_{n=1}^{\infty}\tilde{I}_{A_n}f(x_n)$.

\noindent $(2)$. Let $f:E\rightarrow \bar{L}^{0}(\mathcal{F})$ have the local property and $G\subset E$ be a nonempty subset, then $\bigvee\{f(x)~|~x\in G\}=\bigvee\{f(x)~|~x\in H_{cc}(G)\}$.

\noindent $(3)$. Let $f$ and $g$ be any two functions from $E$ to $\bar{L}^{0}(\mathcal{F})$ such that they both have the local property and $G\subset E$ a nonempty subset. If $f(x)=g(x)$ for all $x\in G$, then $f(x)=g(x)$ for all $x\in H_{cc}(G)$.

\noindent $(4)$. Let $\{f_{\alpha},\alpha\in\Gamma\}$ be a family of functions from $E$ to $\bar{L}^{0}(\mathcal{F})$ such that each $f_{\alpha}$ has the locally property, then $f:E\rightarrow \bar{L}^{0}(\mathcal{F})$ defined by $f(x)=\bigvee\{f_{\alpha}(x)~|~\alpha\in\Gamma\}$ for all $x\in E$, also has the local property.
\end{lemma}

\begin{proof} (1). $f(x)=(\Sigma_{n=1}^{\infty}\tilde{I}_{A_n})f(x)=\Sigma_{n=1}^{\infty}\tilde{I}_{A_n}f(x)=\Sigma_{n=1}^{\infty}\tilde{I}_{A_n}f(\tilde{I}_{A_n}x)= \Sigma_{n=1}^{\infty}\tilde{I}_{A_n}f(\tilde{I}_{A_n}x_n)=\Sigma_{n=1}^{\infty}\tilde{I}_{A_n}f(x_n)$.

(2). Let $\xi=\bigvee\{f(x)~|~x\in G\}$ and $\eta=\bigvee\{f(x)~|~x\in H_{cc}(G)\}$, then $\xi\leq\eta$ is clear, it remains to prove $\eta\leq\xi$. For any $x\in H_{cc}(G)$, let $\Sigma_{n=1}^{\infty}\tilde{I}_{A_n}g_n$ be a canonical representation of $x$, then $f(x)=\Sigma_{n=1}^{\infty}\tilde{I}_{A_n}f(g_n)\leq\xi$, so $\eta\leq\xi$.

(3). It is clear by (1).

(4). It is also clear by definition.
\end{proof}

We can now prove Theorem 5.5.

\noindent{\bf Proof of Theorem 5.5.}
We first consider the special case when ${\mathcal P}$ has the countable concatenation property. Since $E$ has the countable concatenation property, $f$ is also $\mathcal{T}_{\varepsilon,\lambda}$--lower semicontinuous by Lemma 5.6. Further, since $E^{\ast}_{\varepsilon,\lambda}=E^{\ast}_c$ by Corollary 3.6 and it is obvious that $f^{\ast\ast}_{c}=f^{\ast\ast}_{\varepsilon,\lambda}$, the proof follows from Theorem 5.3.

Now, we consider the general case, namely ${\mathcal P}$ may not necessarily have the countable concatenation property. We consider the problem in $(E, {\mathcal P}_{cc})$. Since ${\mathcal P}_{cc}$ has the countable concatenation property and the locally $L^0$--convex topology induced by ${\mathcal P}_{cc}$ is stronger than that induced by ${\mathcal P}$, applying the special case which has been proved above to $f$ and $(E, {\mathcal P}_{cc})$ we can obtain:$$f(x)=\bigvee\{u(x)-f^\ast_c(u)~|~u\in (E, {\mathcal P}_{cc})^\ast_c\}~\hbox{ for all $x\in E$.}$$
Since $f^\ast_c$ has the local property and $u(x)$ is, of course, local with respect to $u$ for a fixed $x\in E$, then $u(x)-f^{\ast}_{c}(u)$ is local with respect to $u$ when $x$ is fixed. So by (2) of Lemma 5.7 and the fact that $(E, {\mathcal P}_{cc})^\ast_c=H_{cc}(E^\ast_c)$ (namely, Theorem 3.7 , where $E^\ast_c=(E, {\mathcal P})^\ast_c$) we have that $f(x)=\bigvee\{u(x)-f^{\ast}_{c}(u)~|~u\in H_{cc}(E^{\ast}_c)\}$  $=\bigvee\{u(x)-f^{\ast}_{c}(u)~|~u\in E^{\ast}_c\}$.
\hfill \hfill $\square$

\subsection{The Fenchel-Moreau dual representation theorems for nonproper functions in random locally convex modules}

In classical convex analysis, people very often need to consider the Fenchel-Moreau dual representation theorem for a not necessarily proper extended real-valued function, where the notion of a closed function is important. Let $(E,\mathcal{T})$ be a locally convex space. $f:E\rightarrow[-\infty,+\infty]$ is closed if either $f\equiv+\infty$, or $f\equiv-\infty$, or $f$ is a proper lower semicontinuous, cf. \cite{ET}. Thus we should also define and study closed functions in the random setting. In fact, D. Filipovi\'{c}, M. Kupper and N. Vogepoth already studied the problem for a special class of $RN$ module $L^{p}_{\mathcal{F}}(\mathcal{E})$ for financial applications. Here, we make use of Theorem 5.3 to give a unified treatment for the problem.

Let $(E,\mathcal{P})$ be an $RLC$ module over $R$ with base $(\Omega,\mathcal{F},P)$ and $f:E\rightarrow \bar{L}^{0}(\mathcal{F})$ an local function. Let us first give the following notation:

$\mathscr{A}=\{A\in \mathcal{F}~|~$there is $x\in E$ such that $\tilde{I}_A f(x)=\tilde{I}_A(-\infty)\}$;

$\mathscr{B}=\{A\in \mathcal{F}~|~\tilde{I}_A f=\tilde{I}_A(+\infty),$ namely $\tilde{I}_A f(x)=\tilde{I}_A(+\infty)$ for all $x\in E\}$;

$MI(f)=esssup(\mathscr{A})$;

$PI(f)=esssup(\mathscr{B})$;

$BP(f)=\Omega\setminus(MI(f)\bigcup PI(f))$;

$\mathscr{D}=\{A\subset BP(f)~|~A\in \mathcal{F}$ is such that there are $D\in\mathcal{F}$ with $D\subset A$ and $x\in E$ satisfying $f(x)<+\infty$ on $D\}$.

Here, $esssup(\mathcal{H})$ denotes the essential supremum of a subfamily $\mathcal{H}$ of $\mathcal{F}$, cf. \cite{FKV,TXG-JFA}. We can think that $MI(f)$ and $PI(f)$ are disjoint.

It is obvious that $\tilde{I}_{PI(f)}f=\tilde{I}_{PI(f)}(+\infty)$ and $f(x)>-\infty$ on $BP(f)$ for all $x\in E$. Since $\mathscr{A}$ and $\mathscr{D}$ are upward directed, one can use the essential supremum theorem to prove Proposition 5.8 below.

\begin{proposition} We have the following statements:

\noindent $(1)$. There are a countable partition $\{A_n,n\in N\}$ of $MI(f)$ to $\mathcal{F}$ and a sequence $\{y_n,n\in N\}$ in $E$ such that $\tilde{I}_{A_n}f(y_n)=\tilde{I}_{A_n}(-\infty)$ for each $n\in N$.

\noindent $(2)$. There are a countable partition $\{D_n,n\in N\}$ of $BP(f)$ to $\mathcal{F}$ and a sequence $\{x_n,n\in N\}$ in $E$ such that $f(x_n)<+\infty$ on $D_n$ for each $n\in N$ (namely, each $\tilde{I}_{D_n}f$ is proper). Further, if, in addition, $P(BP(f))>0$, then each $D_n$ can be chosen such that $P(D_n)>0$.
\end{proposition}

Let us observe that if $E$ has the countable concatenation property then the local property of $f$ can be used to prove: there are $y\in E$ such that $\tilde{I}_{MI(f)}f(y)=\tilde{I}_{MI(f)}(-\infty)$, and $x\in E$ such that $f(x)<+\infty$ on $BP(f)$, namely $\tilde{I}_{BP(f)}f$ is proper.

For each $D\in\mathcal{F}$, let $E^{D}=\tilde{I}_{D}E:=\{\tilde{I}_{D}x~|~x\in E\}$ and $\|\cdot\|^{D}=$ the restriction of $\|\cdot\|$ to $E^{D}$ for each $\|\cdot\|\in\mathcal{P}$. Then $(E^{D},\mathcal{P}^{D})$ can , of course, be regarded as an $RLC$ module over $R$ with base $(D,D\bigcap\mathcal{F},P(\cdot|D))$ if $P(D)>0$, where $\mathcal{P}^{D}=\{\|\cdot\|^{D}~|~\|\cdot\|\in\mathcal{P}\}$. Further, $f_D:E^{D}\rightarrow \tilde{I}_{D}\bar{L}^{0}(\mathcal{F})$ is defined by $f_{D}(\tilde{I}_{D}x)=\tilde{I}_{D}f(\tilde{I}_{D} x)$ for all $x\in E$.

We can now introduce the notion of a closed function. We can assume, without loss of generality, that $P(BP(f))>0$ for the function $f$ in discussion.

\begin{definition} Let $(E,\mathcal{P})$ be an $RLC$ module over $R$ with base $(\Omega,\mathcal{F},P)$, $f:E\rightarrow \bar{L}^{0}(\mathcal{F})$ a local function. Then $f$ is $\mathcal{T}_{\varepsilon,\lambda}$ (resp., $\mathcal{T}_{c}$)-closed if $\tilde{I}_{MI(f)}f=\tilde{I}_{MI(f)}(-\infty)$ and if $f_A$ is
$L^{0}(A\cap\mathcal{F})-$convex and $\mathcal{T}_{\varepsilon,\lambda}$ (resp., $\mathcal{T}_{c}$)-lower semicontinuous for each $A\in\mathcal{F}$ with $A\subset BP(f)$ and $P(A)>0$ such that $f_A$ is proper.
\end{definition}

\begin{remark} First, $A$ in Definition 5.9 universally exists, for example, let $\{D_n,n\in N\}$ be the same as in (2) of Proposition 5.8, then each $f_{D_{n}}$ is proper. Furthermore, if $f$ is a closed function then $f=\tilde{I}_{PI(f)}(+\infty)+\tilde{I}_{MI(f)}(-\infty)+\sum_{n=1}^{\infty}\tilde{I}_{D_n}f$ with each $\tilde{I}_{D_n}f$ (namely $f_{D_n}$) is proper $L^{0}-$convex lower semicontinuous, so our definition of a closed function is not only very similar to the classical definition of a closed function but also more complicated than the latter. By the way, it is easy to see that a closed function must be $L^{0}-$convex. Secondly, the notion of a $\mathcal{T}_c-$closed function in the sense of Definition 5.9 is more general than that introduced in \cite{FKV-appro}: \cite{FKV-appro} only considered the special case when $E=L^{p}_{\cal F}(\cal E)$, in which case $\tilde{I}_{BP(f)}f$ is proper, whereas $\tilde{I}_{BP(f)}f$ is not necessarily proper in our general case and the study of our general case needs a decomposition of $BP(f)$ as in (2) of Proposition 5.8. Besides, \cite{FKV-appro} employed the strongest notion of a $\mathcal{T}_c-$lower semicontinuous function, whereas we employ the weakest one.
\end{remark}

\begin{proposition} Let $(E,\mathcal{P})$ be the same as in Definition 5.9, $\{f_\alpha,\alpha\in\Gamma\}$ a family of $\mathcal{T}_{\varepsilon,\lambda}$ (resp., $\mathcal{T}_{c}$)-closed functions from $(E,\mathcal{P})$ to $\bar{L}^{0}(\mathcal{F})$ and define $f=\bigvee\{f_\alpha:\alpha\in \Gamma\}$ by $f(x)=\bigvee\{f_\alpha(x):\alpha\in \Gamma\}$ for all $x\in E$. Then $f$ is still $\mathcal{T}_{\varepsilon,\lambda}$ (resp., $\mathcal{T}_{c}$)-closed.
\end{proposition}

\begin{proof} It is easy to see that $MI(f)=essinf\{MI(f_{\alpha}),\alpha\in\Gamma\}$, $PI(f)=esssup\{PI(f_{\alpha}),\alpha\in\Gamma\}$ and $\tilde{I}_{MI(f)}f=\tilde{I}_{MI(f)}(-\infty)$. It remains to show that $f_A$ is $L^{0}(A\cap\cal F)-$convex and $\mathcal{T}_{\varepsilon,\lambda}$ (resp. $\mathcal{T}_c$)$-$lower semicontinuous for each $A\in\cal F$ with $A\subset BP(f)$ and $P(A)>0$ such that $f_A$ is proper. We only gives the proof for the $(\varepsilon,\lambda)-$topology since the case for the locally $L^{0}-$convex topology is similar.

Since each $f_{\alpha}$ is $\mathcal{T}_{\varepsilon,\lambda}-$closed, each $f_{\alpha}$ is $L^{0}-$convex, then $f$ is $L^{0}-$convex, so $f_A$ is $L^{0}(A\cap\cal F)-$convex. Further, since $epi(f_A)=\cap_{\alpha\in\Gamma}epi((f_{\alpha})_A)$, we only need to check that each $epi((f_{\alpha})_A)$ is $\mathcal{T}_{\varepsilon,\lambda}-$closed in $\tilde{I}_A(E\times L^{0}(\cal F))$. In fact, for any fixed $\alpha\in\Gamma$, $A$ must be a subset of $(PI(f_{\alpha}))^c$ since $A\subset BP(f)$, so $A=(A\cap BP(f_{\alpha}))\cup(A\cap MI(f_{\alpha}))$. According to the fact that $\tilde{I}_{MI(f_{\alpha})}f_{\alpha}=\tilde{I}_{MI(f_{\alpha})}(-\infty), epi((f_{\alpha})_A)=epi((f_{\alpha})_{A\cap BP(f_{\alpha})})+\tilde{I}_{A\cap MI(f_{\alpha})}(E\times L^{0}(\cal F))$. Since $f_A$ is proper, it is obvious that $(f_{\alpha})_{A\cap BP(f_{\alpha})}$ is also proper, which shows that $epi((f_{\alpha})_{A\cap BP(f_{\alpha})})$ is $\mathcal{T}_{\varepsilon,\lambda}-$closed in $\tilde{I}_{A\cap BP(f_{\alpha})}(E\times L^{0}(\cal F))$ since $f_{\alpha}$ is a $\mathcal{T}_{\varepsilon,\lambda}-$closed function. Again by noting the fact that $A\cap BP(f_{\alpha})$ and $A\cap MI(f_{\alpha})$ are disjoint we have that $epi((f_{\alpha})_A)$ is $\mathcal{T}_{\varepsilon,\lambda}-$closed.
\end{proof}

\begin{definition} Let $(E,\mathcal{P})$ and $f$ be the same as in Definition 5.9. The greatest $\mathcal{T}_{\varepsilon,\lambda}$ (resp., $\mathcal{T}_{c}$)-closed function majorized by $f$, denoted by $Cl_{\varepsilon,\lambda}(f)$ (resp., $Cl_c(f)$), is the $\mathcal{T}_{\varepsilon,\lambda}$ (resp., $\mathcal{T}_c$)-closure of $f$.
\end{definition}

\begin{lemma} Let $(E,\mathcal{P})$ and $f$ be the same as in Definition 5.9. If $f$ is $\mathcal{T}_{\varepsilon,\lambda}$--closed, then $f^{\ast\ast}_{\varepsilon,\lambda}=f$.
\end{lemma}

\begin{proof} Since $f$ is $\mathcal{T}_{\varepsilon,\lambda}$--closed, it is obvious that $\tilde{I}_{MI(f)}f^{\ast\ast}_{\varepsilon,\lambda}=\tilde{I}_{MI(f)}f=\tilde{I}_{MI(f)}(-\infty)$ and $\tilde{I}_{PI(f)}f^{\ast\ast}_{\varepsilon,\lambda}=\tilde{I}_{PI(f)}f=\tilde{I}_{PI(f)}(+\infty)$. Let $\{D_n,n\in N\}$ be the same as in (2) of Proposition 5.8 with $P(D_n)>0$ for all $n\in N$,then each $f_{D_n}$ is a proper $L^{0}(D_n\bigcap\mathcal{F})$--convex $\mathcal{T}_{\varepsilon,\lambda}$--lower semicontinuous on $E^{D_n}$. It is also obvious that $\tilde{I}_{D_n}f^{\ast\ast}_{\varepsilon,\lambda}=f^{\ast\ast}_{D_n}=f_{D_n}=\tilde{I}_{D_n} f$ for each $n\in N$ by Theorem 5.3, so $f^{\ast\ast}_{\varepsilon,\lambda}=f$.
\end{proof}

\begin{theorem} Let $(E,\mathcal{P})$ and $f$ be the same as in Definition 5.9. Then $f^{\ast\ast}_{\varepsilon,\lambda}=Cl_{\varepsilon,\lambda}(f)$.
\end{theorem}

\begin{proof} It is obvious that $f^{\ast\ast}_{\varepsilon,\lambda}\leq f$ and $f^{\ast\ast}_{\varepsilon,\lambda}$ is $\mathcal{T}_{\varepsilon,\lambda}$--closed, so $f^{\ast\ast}_{\varepsilon,\lambda}\leq Cl_{\varepsilon,\lambda}(f)$. On the other hand, $Cl_{\varepsilon,\lambda}(f)\leq f$, then $Cl_{\varepsilon,\lambda}(f)=(Cl_{\varepsilon,\lambda}(f))^{\ast\ast}_{\varepsilon,\lambda}\leq f^{\ast\ast}_{\varepsilon,\lambda}$ by Lemma 5.13.
\end{proof}

\begin{corollary} Let $(E,\mathcal{P})$ be an $RLC$ module over $R$ with base $(\Omega,\mathcal{F},P)$ such that $E$ has the countable concatenation property and $f:E\rightarrow\bar{L}^{0}(\mathcal{F})$ a local function, then $f^{\ast\ast}_{c}=Cl_{c}(f)$.
\end{corollary}

\begin{proof} It is similar to the proof of Theorem 5.5, so is omitted. \end{proof}

\sec{Acknowledgements:} {The first author of this paper thanks Professor Quanhua Xu for kindly providing us the excellent literature \cite{HLR} in February, 2012 to make us know, for the first time, the existence of \cite{HLR} since \cite{HLR} has never been mentioned before in the literature of related fields (e.g. \cite{SS}).}


\end{document}